\documentclass{svmult}

\usepackage{makeidx}         
\usepackage{graphicx}        
\usepackage{multicol}        
\usepackage[bottom]{footmisc}

\makeindex             

\usepackage{amssymb}                   
\usepackage[intlimits]{amsmath}  

\newcommand{\req}[1]{(\ref{eq:#1})}

\renewcommand{\S}{\mathcal{S}}
\renewcommand{\D}{\mathcal{D}}

\newcommand\R{{\mathbb{R}}}
\newcommand\N{{\mathbb{N}}}
\newcommand\C{{\mathbb{C}}}

\newcommand{\ve}{\varepsilon}

\newcommand{\vx}{{\vec{x}}}

\renewcommand{\k}{\mathbf{k}}
\newcommand\vO{{\mathbf{0}}}
\newcommand\ta{\tilde{\alpha}}

\newcommand{\patt}{p_{\text{att}}}
\newcommand{\sgn}{\mbox{sgn}}
\newcommand{\Heavi}{H}
\newcommand{\ef}[1]{\exp \left( #1 \right)}
\newcommand{\M}{\mathcal{M}}

\renewcommand{\i}{{\mathrm{i}}\,}

\newcommand{\supp}{\mbox{supp}}
\newcommand{\F}{\mathcal{F}}
\renewcommand{\H}{\mathcal{H}}
\newcommand{\A}{\mathcal{A}}
\newcommand{\B}{\mathcal{B}}

\newcommand{\set}[1]{\{ #1 \}}

\renewcommand{\d}{\mbox{d}}
\def\Xint#1{\mathchoice
  {\XXint\displaystyle\textstyle{#1}}%
  {\XXint\textstyle\scriptstyle{#1}}%
  {\XXint\scriptstyle\scriptscriptstyle{#1}}%
  {\XXint\scriptscriptstyle\scriptscriptstyle{#1}}%
  \!\int}
\def\XXint#1#2#3{{\setbox0=\hbox{$#1{#2#3}{\int}$}
    \vcenter{\hbox{$#2#3$}}\kern-.5\wd0}}

\newcommand{\abs}[1]{\left|#1\right|}

\newcommand{\hilbert}[1]{\mathcal{H} \left\{ #1 \right\}}
\newcommand{\ilaplace}[1]{\mathcal{L}^{-1} \left\{ #1 \right\}}

\newcommand{\fourier}[1]{\mathcal{F} \left\{ #1 \right\}}
\newcommand{\ifourier}[1]{\mathcal{F}^{-1} \left\{ #1 \right\}}

\newcommand{\green}{\mathcal{G}}

\newcommand{\ltext}[1]{\qquad \text{#1} \qquad}
\newcommand{\stext}[1]{\quad \text{#1} \quad}

\RequirePackage{ifthen}
\makeatletter
\newcommand{\logmessage}[1]{\@latex@warning{#1}}
\newcommand{\ignore}{\logmessage{Text ignored}\@gobble} \makeatother

  \newtheorem{notation}[theorem]{Notation \& Terminology}

\begin{document}

\title*{Photoacoustic Imaging Taking into Account Attenuation}
\titlerunning{Photoacoustics}
\author{Richard Kowar\inst{1}\and
Otmar Scherzer \inst{2}}
\institute{Institute of Mathematics, University of Innsbruck, Technikerstr. 21a, A-6020 Innsbruck, Austria
\texttt{richard.kowar@uibk.ac.at} \\
Computational Science Center, University of Vienna, Nordbergstr. 15, A-1090 Vienna, Austria
\texttt{otmar.scherzer@univie.ac.at}}

\maketitle

\section{Introduction}

\emph{Photoacoustic Imaging} is one of the recent hybrid imaging techniques,
which attempts to visualize the distribution of the \emph{electromagnetic absorption coefficient}
inside a biological object. In photoacoustic experiments, the medium is exposed to a short pulse
of a relatively low frequency electromagnetic (EM) wave. The exposed medium absorbs a fraction of the
EM energy, heats up, and reacts with thermoelastic expansion. This induces acoustic waves,
which can be recorded outside the object and used to determine the electromagnetic absorption coefficient.
The combination of EM and ultrasound waves (which explains the usage of the term \emph{hybrid}) allows one to
combine high contrast in the EM absorption coefficient with high resolution of ultrasound.
The method has demonstrated great potential for biomedical applications, including functional brain imaging of animals
\cite{WanPanKuXieSto03}, soft-tissue characterization, and early stage cancer diagnostics \cite{KruKisMilRey00},
as well as imaging of vasculature \cite{ZhaLauBea07}. For a general survey on biomedical applications see \cite{XuWan06}.
In comparison with the X-Ray CT, photoacoustics is non-ionizing. Its further advantage
is that soft biological tissues display high contrasts in their ability to absorb frequency electromagnetic waves.
For instance, for radiation in the near infrared domain, as produced by a Nd:YAG laser, the absorption coefficient in
human soft tissues varies in the range of $0.1/\textrm{cm}$--$0.5/\textrm{cm}$ \cite{ChePraWel90}.
The contrast is also known to be high between healthy and cancerous cells, which makes photoacoustics a promising early
cancer detection technique. Another application arises in biology: Multispectral optoacoustic tomography
technique is capable of high-resolution visualization of fluorescent proteins deep within highly light-scattering living organisms
\cite{RazDisVinMaPerKoeNtz09}. In contrast, the current fluorescence microscopy techniques are limited to the depth of
several hundred micrometers, due to intense light scattering.

Different terms are often used to indicate different excitation sources: \emph{Optoacoustics} refers to illumination in the
visible light spectrum, \emph{Photoacoustics} is associated with excitations in the visible and infrared range, and \emph{Thermoacoustics}
corresponds to excitations in the microwave or radio-frequency range. In fact, the carrier frequency of the illuminating pulse is varying,
which is usually not taken into account in mathematical modeling. Since the corresponding mathematical models are equivalent,
in the mathematics literature, the terms opto-, photo-, and thermoacoustics are used interchangeably.
In this article, we are addressing only the \emph{photoacoustic tomographic technique} PAT (which is mathematically
equivalent to the thermoacoustic tomography TAT).

Various kinds of photoacoustic imaging techniques have been implemented. One should distinguish between photoacoustic
\emph{microscopy} (PAM) and \emph{tomography} (PAT). In microscopy, the object is scanned pixel by pixel (or voxel by voxel).
The measured pressure data provides an image of the electromagnetic absorption coefficient \cite{ZhaMasStoWan06}.
Tomography, on the other hand, measures pressure waves with detectors surrounding completely or partially the object. Then the
internal distribution of the absorption coefficients is reconstructed using mathematical inversion techniques (see the sections below).

The common underlying mathematical equation of PAT is the \emph{wave equation} \index{wave equation, photoacoustic model}
for the pressure
\begin{equation}\label{eq:ex:wave3d}
\boxed{
\frac{1}{c_0^2} \frac{\partial^2 p}{\partial t^2}(\vx,t ) - \nabla^2 p(\vx,t) = \frac{d j}{d  t}(t)
      \left(  \frac{\mu_{\rm abs}(\vx) \beta (\vx) J({\vx})}{c_p(\vx)} \right)\,,\;\vx \in \R^3,\,t > 0 \;.}
\end{equation}
Here $c_p$ denotes the specific heat capacity, $J$ is the spatial intensity distribution, $\mu_{\rm abs}$
denotes the absorption coefficient, $\beta$ denotes the thermal
expansion coefficient and $c_0$ denotes the speed of sound, which is commonly assumed to be constant. The
assumption that there is no acoustic pressure before the object is
illuminated at time $t  = 0$ is expressed by
\begin{equation}\label{eq:ex:ini3d}
\boxed{
    p(\vx ,t) = 0 \,, \qquad  \vx \in \R^3, t < 0\;.}
\end{equation}
In PAT, $j(t)$ approximates a pulse, and can be considered as a $\delta$-impulse $\delta(t)$. Introducing the
shorthand notations \index{absorption energy}
\begin{equation}\label{eq:u(x)}
\begin{aligned}
\rho(\vx) := \frac{\mu_{\rm abs}(\vx) \beta (\vx) J(\vx)}{c_p(\vx)}\,,
\end{aligned}
\end{equation}
one reduces (\ref{eq:ex:wave3d}) and (\ref{eq:ex:ini3d}) to
\begin{equation}\label{eq:ex:ivp3d}
\boxed{
\frac{1}{c_0^2} \frac{\partial^2 p}{\partial t^2} (\vx,t) - \nabla^2 p(\vx,t)= 0\,,\quad  \vx \in \R^3,  t  > 0\,,}
\end{equation}
\index{wave equation}
with initial values
\begin{equation}
\label{eq:init_values}
\boxed{
   p(\vx,0) = \rho(\vx)\,, \quad \frac{\partial p}{\partial t}(\vx,0) = 0 \qquad  \vx \in \R^3\;.}
\end{equation}
The quantity $\rho$ in (\ref{eq:ex:wave3d}) and (\ref{eq:u(x)}) is a combination of several physical parameters. All
along this paper $\rho$ should not be confused with the source term
\begin{equation}
\label{eq:source}
\boxed{f(\vx,t) = \frac{d j}{d  t}(t) \rho(\vx)\,,\quad  \vx \in \R^3,  t  > 0  \;.}
\end{equation}
\index{source term, photoacoustic model}
In PAT, some data about the pressure $p(\vx,t)$ are measured and the main task is to reconstruct the initial pressure $\rho$ from these data.
While the excitation principle is always as described above and thus (\ref{eq:ex:ivp3d}) holds, the specific type of data measured
depends on the type of transducers used, and thus influences the mathematical model.

Nowadays there is a trend to incorporate more and more modeling into photoacoustic.
In particular, taking into account locally varying \emph{wave speed} and \emph{attenuation}.
Even more there is a novel trend to \emph{qualitative photoacoustics}, which is concerned with
estimating physical parameters from the imaging parameter of standard photoacoustics.
In this paper we focus on attenuation correction, where we survey some recent
progress. Inversion with varying wave speed has been considered for instance in
\cite{AgrKuc07,HriKucNgu08}, and is not further discussed here.

The outline of this paper is as follows: First, we review existing attenuation models and discuss their causality properties, which we
believe to be essential for algorithms for inversion with attenuated data. Then, we survey causality properties of common attenuation models.
We also derive integro-differential equations which the attenuated waves are satisfying. In addition we discuss the ill--conditionness of the
inverse problem for calculating the unattenuated wave from the attenuated one.

\section{Attenuation}
\label{sec:attenuation}
The difficult issue of effects of and corrections for the attenuation of acoustic waves
in PAT has been studied\cite{RivZhaAna06,BurGruHalNusPal07,PatGre06,KowSchBon10},
although no complete conclusion on the feasibility of these models has been reached.

Mathematical models for describing attenuation are formulated in the
frequency domain, taking into account that attenuation disperses
high frequency components more rapidly over traveled distance. Let
$\green(\vx,t)$ denote the attenuated wave which
originates from an impulse ($\delta_{\vx,t}$-distribution) at $\vx=0$ at
time $t=0$. In mathematical terms $\green$ is the Green-function of
attenuated wave equation. Moreover, we denote by
\begin{equation}
\boxed{
\label{eq:g0} \green_0(\vx,t) =
\frac{\delta\left(t-\frac{\abs{\vx}}{c_0}\right)}{4\,\pi\,\abs{\vx}}}\index{Green function, standard wave equation}
\end{equation}
the Green function of the unattenuated wave equation;
that is, it is the solution of \eqref{eq:ex:ivp3d},
\eqref{eq:init_values} with constant sound speed $c(x) \equiv c_0$
and initial conditions
$$
   \green_0(\vx,0)=0 \quad\mbox{ and}\quad
   \frac{\partial \green_0}{\partial t}(\vx,0)=\delta_{\vx,t}\,.
$$
Common mathematical formulations of \emph{attenuation} assume that \index{attenuation}
\index{attenuation coefficient, $\beta^*$}
\begin{equation}
\label{eq:GG0}
\boxed{
\fourier{\green}(\vx,\omega) =  \ef{ -\beta^*
(\abs{\vx},\omega)}\, \fourier{\green_0}(\vx,\omega)\,, \quad \vx
\in \R^3,\, \omega \in \R\;.}
\end{equation}
Here $\fourier{\cdot}$ denotes the Fourier transform with respect to time $t$ (cf. Appendix~\ref{sec:app}).
Applying the inverse Fourier transform $\ifourier{\cdot}$ to \req{GG0} gives
\begin{equation}
\label{eq:FGG0}
\boxed{\green(\vx,t) =   K(\vx,t) *_t \green_0(\vx,t)}
\qquad \mbox{($*_t$ time convolution)}
\end{equation}
\index{convolution kernel, $K$}
where
\begin{equation}
\label{eq:kernel}
K(\vx,t) := \frac{1}{\sqrt{2\,\pi}} \ifourier{\ef{-\beta^* (\abs{\vx},\cdot)}}(t)\;.
\end{equation}
From (\ref{eq:FGG0}) and (\ref{eq:g0}) it follows that
\begin{equation*}
\begin{aligned}
\green (\vx,t) &= K(\vx,t) *_t \green_0(\vx,t)\\
&= \int_\R K(\vx,t-\tau) \frac{\delta (\tau - \frac{\abs{\vx}}{c_0})}{4 \pi \abs{\vx}}\,d\tau \\
&= \frac{K\left(\vx,t-\frac{\abs{\vx}}{c_0}\right)}{4\pi \abs{\vx}}\;.
\end{aligned}
\end{equation*}
Consequently,
\begin{equation}
\label{eq:Kandgreen}
\boxed{
\green (\vx,t+\abs{\vx}/c_0) = K(\vx,t)/(4\pi \abs{\vx})\;.}
\end{equation}
Moreover, we emphasize that the Fourier transform of a real and even (real and
odd) function is real and even (imaginary and odd). Since $\green$
and $\green_0$ are real valued, $K$ must be real valued and
consequently the real part $\Re(\beta^*)$ of $\beta^*$ has to be
even with respect to the frequency $\omega$ and $\Im(\beta^*)$ has
to be odd with respect to $\omega$. Attenuation is caused if
$\Re(\beta^*)$ is positive and since then $\beta^*$ has a nonzero
imaginary part due to the Kramers-Kronig relation, attenuation causes dispersion.
In the literature the
following product ansatz is commonly used
\begin{equation}
\label{eq:attenuation_law}
\boxed{
\beta^* (\abs{\vx},\omega) = \alpha^*(\omega)\,\abs{\vx} \qquad \omega\in\R,\,\vx\in\R^3\;.}
\end{equation}
In the sequel we concentrate on these models and use the following
terminology:
\begin{definition}
We call  $\beta^*$ of standard form if~(\ref{eq:attenuation_law}) holds. Then the function
\index{attenuation coefficient, standard}
\begin{equation}
\label{eq:attenuation_coefficient}
\alpha^*: \R \to \C
\end{equation}
is called \emph{standard attenuation coefficient} and $\alpha = \Re(\alpha^*)$ is
called the \emph{attenuation law}. \index{attenuation law} We also call $\beta^*$ the
\emph{attenuation coefficient}. \index{attenuation coefficient}
\end{definition}
From the relation (\ref{eq:attenuation_law}), it follows that $\Re (\alpha^*)$ is even,
$\Im (\alpha^*)$ is odd, and $\Re (\alpha^*)>0$ (the last inequality guarantees attenuation).

In the following we summarize common attenuation coefficients and
laws:
In what follows $\alpha_0$ denotes a positive parameter and
\begin{equation}
\label{ta0} \ta_0 = \frac{\alpha_0}{\cos
\left(\frac{\pi}{2}\gamma\right)}           \qquad\quad (0 < \gamma \not\in \N)\,
\end{equation}
is a possibly non-positive coefficient.
\begin{itemize}
\item {\bf Frequency Power Laws:}
\begin{itemize}
\item Let $0 < \gamma \not\in \N$.
The frequency power law \emph{attenuation coefficient} is defined by
\begin{equation}
\label{eq:powlaw1}
\alpha_{pl}^*(\omega) = \ta_0 (-\i \omega)^\gamma = \ta_0 \abs{\omega}^\gamma \left(
\cos \left(\frac{\pi}{2} \gamma \right) - \i  \sgn(\omega) \sin \left(\frac{\pi}{2} \gamma \right)\right)
\end{equation}
for $\omega \in \R$.
Therefore, the \emph{attenuation law} is given by \index{attenuation law, power law, $\alpha_{pl}^*$}
\begin{equation}
\label{eq:powlaw1b}
\boxed{
\alpha_{pl}(\omega) = \alpha_0\,\abs{\omega}^\gamma\;.}
\end{equation}
These models have been considered for instance in \cite{Sza94,Sza95,WatHugBraMil00,WaMoMi05}.
\item Let $\gamma=1$ and $\omega_0\neq0$, the attenuation coefficient is defined by
\begin{equation}\label{alphagamma1}
\begin{aligned}
  \alpha_{pl}^*(\omega)
     := \alpha_0\,\abs{\omega}
        + \i \frac{2}{\pi}\,\alpha_0\,\omega\,\log\left|\frac{\omega}{\omega_0}\right|
\qquad \omega\in\R\,.
\end{aligned}
\end{equation}
The attenuation law is
\index{attenuation law, power law, $\alpha_{pl}^*$}
\begin{equation}
\label{alphagamma1a}
\boxed{
\alpha_{pl}(\omega):=\alpha_0\,\abs{\omega}\;.}
\end{equation}
This model has been considered in \cite{Sza95,WaMoMi05}.
\end{itemize}

\item{{\bf Szabo:}} Let $0 < \gamma\not\in\N$. The attenuation coefficient
\footnote{In this paper the root of a complex number is always the one with
         non-negative real part.} of
Szabo's law is defined by
\begin{equation}\label{alpha*szabo}
\begin{aligned}
  \alpha_{sz}^*(\omega) &= \frac{1}{c_0}\,\sqrt{ (-\i \omega)^2
        + 2 \ta_0 c_0 (-\i \omega)^{\gamma+1}  }
         + \i \frac{\omega}{c_0}\;.
\end{aligned}
\end{equation}
We denote Szabo's attenuation law by \index{attenuation law, Szabo, $\alpha_{sz}^*$}
\begin{equation*}
\boxed{
\alpha_{sz}(\omega) := \Re (\alpha_{sz}^*(\omega))\;.}
\end{equation*}
For small frequencies $\alpha_{sz}(\omega)$ behaves like
$\alpha_0\,|\omega|^\gamma$.
This model has been considered in \cite{Sza94,Sza95} where, in addition, also a model for $\gamma \in \N$
has been introduced.

\item{{\bf Thermo-Viscous Attenuation Law:} (see e.g.~\cite{KinFreCopSan00,Sza94}):} Here, for
$\tau_0>0$, the attenuation coefficient is defined by
\begin{equation}\label{alpha*th}
\begin{aligned}
  \alpha_{tv}^*(\omega)
     =  \frac{-\i \omega}{c_0\,\sqrt{1-\i \tau_0\,\omega}}
         +  \frac{\i \omega}{c_0}\;
\end{aligned}
\end{equation}
with attenuation law \index{attenuation law, thermo-viscous, $\alpha_{tv}^*$}
\begin{equation}\label{alpha*threal}
\begin{aligned}
\boxed{
  \alpha_{tv}(\omega)
     =  \frac{\tau_0\, \omega^2}{ \sqrt{2}\,c_0\,\sqrt{ (1+\sqrt{1+(\tau_0\,\omega)^2})\,(1+(\tau_0\,\omega)^2)} } \;.}
\end{aligned}
\end{equation}
For small frequencies $\alpha_{tv}(\omega)$ behaves like
$\frac{\tau_0\,\omega^2}{2\,c_0}$. That is the thermo-viscous law
approximates a power attenuation law with exponent $2$.

\item {{\bf Nachman, Smith and Waag \cite{NacSmiWaa90}:}}
Consider a homogeneous and isotropic fluid with density $\rho_0$ in which $N$ relaxation processes take place.
Then the attenuation coefficient of the model in  \cite{NacSmiWaa90} reads as follows:
\begin{equation}\label{alpha*Nachman+}
\begin{aligned}
  \alpha_{nsw}^*(\omega)
     =  \frac{-\i \omega}{c_0}\,\left[
             \frac{c_0}{\tilde c_0}\, \sqrt{ \frac{1}{N}\, \sum_{m=1}^N
          \frac{1-\i\,\tilde\tau_m\,\omega}{1-\i \tau_m\,\omega} }
            -1 \right]\;.
\end{aligned}
\end{equation}
All parameters appearing in (\ref{alpha*Nachman+}) are positive and real.
$\kappa_m$ and $\tau_m$ denote the compression modulus and the relaxation
time of the $m-$th relaxation process, respectively, and
\begin{equation}\label{deftildec0tau0}
  \tilde c_0 := \frac{c_0}{\sqrt{1+\sum_{m=1}^N c_0^2\,\rho_0\,\kappa_m}}
\quad\mbox{and}\quad
  \tilde \tau_m := \tau_m\,(1-N\,\tilde c_0^2\,\rho_0\,\kappa_m)\,.
\end{equation}
The last two definitions imply that
\begin{equation}\label{proptildec0tau0}
  \frac{\tilde c_0^2}{c_0^2} = \frac{1}{N}\,\sum_{m=1}^N \frac{\tilde \tau_m}{\tau_m} \;.
\end{equation}
We denote the according attenuation law by \index{attenuation law, Nachman \& Smith \& Waag, $\alpha_{nsw}^*$}
\footnote{In \cite{NacSmiWaa90} they use the notion $c_\infty$ for $c_0$ and $c$ for $\tilde c_0$.}
\begin{equation*}
\boxed{
\alpha_{nsw}(\omega) := \Re (\alpha_{nsw}^*(\omega))\;.}
\end{equation*}

\item{{\bf Greenleaf and Patch \cite{PatGre06}}} consider for $\gamma \in \set{1,2}$
the attenuation coefficient \index{attenuation law, Greenleaf \& Patch, $\alpha_{gp}^*$}
\begin{equation*}
      \alpha_{gp}^*(\omega) = \alpha_0\,\abs{\omega}^\gamma\,,
\end{equation*}
which, since it is real, equals the attenuation law
\begin{equation}\label{alpha*Patch}
\boxed{
      \alpha_{gp}(\omega) = \Re(\alpha_{gp}^*(\omega))\;.}
\end{equation}
\item{\bf Chen and Holm \cite{CheHolm04}:}
This model describes the attenuation as a function of the absolute value of the vector-valued wave number $\k\in\R^3$
(instead of the frequency $\omega\in\R$). Let $\F_{3D}$ denote the $3D-$Fourier transform
\begin{equation*}
\begin{aligned}
   \F_{3D} \set{ f(\k) }(\vx)
        = \frac{1}{\sqrt{(2\,\pi)^3}}\,\int_{\R^3} \ef{\i \vx\cdot\k}\,  f(\k) \,\d\k \,,
\end{aligned}
\end{equation*}
then the Green function of the attenuated equation is defined by
\begin{equation}\label{GreenChenHolm}
\begin{aligned}
    \green(\vx,t) = \frac{H(t)\,c_0^2}{(2\,\pi)^{3/2}}\,\F_{3D}\left\{ \ef{A(\cdot)\,t}\, \frac{\sin(B(\cdot)\,t)}{B(\cdot)} \right\}(\vx)
\end{aligned}
\end{equation}
where, for given $\alpha_1>0$,
\begin{equation}
\label{AkBk}
    A(\k) := -\alpha_1\,c_0\,\abs{\k}^\gamma\,,
\qquad
    B(\k) := c_0\,\sqrt{\abs{\k}^2-\alpha_1^2\,\abs{\k}^{2\,\gamma}}\,.
\end{equation}

\item{\bf In \cite{KowSchBon10}} we proposed \index{attenuation law, \cite{KowSchBon10}, $\alpha_{ksb}^*$}
\begin{equation}\label{eq:powlaw2}
\begin{aligned}
   \alpha_{ksb}^*(\omega)
       = \frac{\alpha_0\,(-\i \omega)}{c_0\,\sqrt{1+(-\i \tau_0\,\omega)^{\gamma-1}}}
\qquad\qquad (\gamma\in (1,2],\,\tau_0>0)\,,
\end{aligned}
\end{equation}\label{eq:powlaw2a}
where the square root is again the complex root with positive real part.

Let $\gamma\in (1,2]$. Then, for small frequencies we have
\begin{equation*}
\boxed{\alpha_{ksb}(\omega) \approx \frac{\alpha_0\,\sin(\frac{\pi}{2}(\gamma-1))}{2\,c_0\,\tau_0} \,|\tau_0\,\omega|^\gamma >0\,.}
\end{equation*}
Thus our model behaves like a power law for small frequencies.
\end{itemize}

Distinctive features of unattenuated wave propagation (,i.e. the
solution of the standard wave equation) are \emph{causality} and
\emph{finite wave front velocity}. It is reasonable to assume that
the attenuated wave satisfies the same distinctive properties as
well. In the following we analyze causality properties of the
standard attenuation models.

\section{Causality}
In the following we present some abstract definitions and basic
notations. In the remainder $\vx$ will always denote a vector in
three dimensional space. When we speak about functions, we always
mean generalized functions, such as for instance distributions or
tempered distributions - we recall the definitions of (tempered)
distribution in the course of the paper.

\begin{definition}
\label{def:defAc} A function $f:=f(\vx,t)$ defined on the Euclidean
space over time (i.e. in $\R^4$) is said to be \emph{causal} if it
satisfies $f(\vx,t) = 0$ for $t < 0$. \index{Function, causal}
\end{definition}

\begin{notation}
Let $\A: D \to D$ be a linear operator, where $\emptyset \neq D$ is
an appropriate set of functions from $\R^4$ to $\R$. In this paper
we always assume that $\A$ satisfies the following properties:
\begin{itemize}
\item $\A$ is \emph{shift invariant} in space and time. That is, for every function $f$ and every shift
      $L:=L(\vx,t):=(\vx-\vx_0,t-t_0)$, with
      $\vx_0 \in \R^3$ and $t_0 \in \R$, it holds that
      \begin{equation*}
      \A (f \circ L)  = (\A f) \circ L\;.
      \end{equation*}
      \index{Operator, shift invariant}
\item $\A$ is \emph{rotation invariant} in space. That is, for every function $f$ and every
      rotation matrix $R$, it holds that
      \begin{equation*}
      \A (R f)  = R(\A f)\;.
      \end{equation*}\index{Operator, rotation invariant}
\item $\A$ is \emph{causal}. That is, it maps causal functions to causal functions. From
(\ref{eq:a_conv}) it follows that $\A$ is causal, if and only if the
associated Green function is causal. \index{Operator, causal}
\end{itemize}
\begin{definition}
\index{Green function}
\label{def:defA} The \emph{Green function} of $\A$ is defined by
$$\green:=\green(\vx,t)=\A \delta_{\vx,t} (\vx,t)\;.$$
\end{definition}

\begin{remark}
The operator $\A$ is uniquely determined by $\green$ and vice versa.
This follows from the fact that
\begin{equation}
\label{eq:a_conv}
\begin{aligned}
\A f (\vx_0,t_0) &= \A \left( \int_\R \int_{\R^3} f(\vx_0-\vx,t-t_0)
\delta_{\vx,t} (\vx,t)\,d \vx dt \right) \\
&= \int_\R \int_{\R^3} f(\vx_0-\vx,t-t_0) \green(\vx,t)\,d \vx dt\,.
\end{aligned}
\end{equation}
\end{remark}

Moreover, we use the following terminology and abbreviations:
\begin{itemize}
\item From the rotation invariance of $\A$ it follows that
      \begin{equation}
      \label{eq:T}
      \hat{T}(\vx):=\sup \set{t : \green(\vx,\tau) = 0 \text{ for all } \tau \leq t}\,,
      \end{equation}
      is rotationally symmetric, which allows us to use the shorthand notation
      \begin{equation}\label{defT}
      T(r) = \hat{T}(\vx) \text{ where } r = \abs{\vx}\;.
      \end{equation}
      With this notation (\ref{eq:T}) can be equivalently expressed as
      \begin{equation}\label{eq:defT(r)0a}
      \green\left(\vx, t + T(\abs{\vx}) \right) = 0 \text{ for every } t<0\;.
      \end{equation}
      In physical terms $T(\abs{\vx})$ denotes the \emph{travel time of a wave front} originating
      at position $\mathbf{0}$ at time $t=0$ and traveling to $\vx$.
      \index{travel time of a wave front}
\item Because $\green$ is rotationally symmetric we can write
     \begin{equation*}
     \green(\vx,T(\abs{\vx})) = \hat{\green}(r,T(r)) \text{ with } r = \abs{\vx}\;.
     \end{equation*}
     Taking the inverse function of $T$, which we denote by $S=S(t)$, we then find
     \begin{equation*}
     \green(\vx,T(\abs{\vx})) = \hat{\green}(S(t),t)\,,
     \end{equation*}
\item The \emph{wave front} is the set \index{wave front}
      \begin{equation*}
      {\cal W}:=\set{(\vx,T(\abs{\vx})): \vx \in \R^3}\;
      \end{equation*}
\item The \emph{wave front speed} \index{wave front speed} $V$ is the variation of the location of the
      wave front as a function of time. That is,
      \begin{equation}
      \label{eq:v}
      V(t)=\frac{\d S}{\d t}(t)= \left. \frac{1}{T'(r)}\right|_{r=S(t)}\,.
      \end{equation}
      Here $T'$ denotes the derivative with respect to the radial component $r$.
\item We say that $\A$ has a \emph{finite speed of propagation} \index{finite speed of propagation} if there exists a constant
      $\hat{c}_0$ such that
      \begin{equation} \label{eq:defT(r)}
      0 < \left(T'(r)\right)^{-1} \leq  \hat{c}_0 < \infty\;.
      \end{equation}
      In this case it follows from (\ref{eq:v}) that the wave front
      velocity satisfies
      \begin{equation}
      \label{eq:v2}
      V(t) \leq \hat{c}_0 < \infty\;.
      \end{equation}
\item We call an operator $\A$ strongly causal, if it is causal and
satisfies the finite propagation speed property. \index{Operator, strongly causal}
\end{itemize}
\end{notation}

The following lemma addresses the case of attenuation coefficients of standard form and
gives examples of strongly causal operators $\A$.

\begin{lemma}\label{lemm:cconst}
Let  $\beta^*(\abs{\vx},\omega) = \alpha^*(\omega) \abs{\vx}$ be of the standard form~(\ref{eq:attenuation_law}) and $\A$ (\ref{eq:a_conv}) be the
operator defined by the Green function $\green$, which is defined in~(\ref{eq:FGG0}). Then $\A$ is strongly causal if and only if for every
$\vx \in \R^3$ the function
\begin{equation*}
t \to \frac{1}{\sqrt{2\pi}} \ifourier{\ef{- \alpha^*(\omega)\abs{\vx}}}\,,
\end{equation*}
defined in~(\ref{eq:kernel}), is causal.
\end{lemma}

\begin{proof}
We assume that $\A$ is strongly causal. It follows from \cite[Theorem 3.1]{KowSchBon10}
that there exists a constant $c$, which is smaller than or equal to the wave speed $c_0$ from (\ref{eq:ex:ivp3d}),
which satisfies $T(\abs{\vx}) = \frac{\abs{\vx}}{c}$ for all $\vx \in \R^3$.
Using the definitions of the travel time $T(\abs{x})$ and (\ref{eq:kernel}), it follows from (\ref{eq:Kandgreen})
that $t\to K(\vx,t)$ is causal.

Now, for every $\vx \in \R^3$ let $K$ be causal. Then from (\ref{eq:Kandgreen}) it follows
that $t \to \green\left(\vx,t+\abs{\vx}/c_0)\right)$ is causal. Since $T(\abs{\vx})$ denotes the largest positive
time period for which $t \to \green\left(\vx,t+T(\abs{\vx})\right)$ is causal, we have for $r>0$:
\begin{equation}\label{eq:lem1Tc0}
     0 <  \frac{r}{c_0} \leq   T(r) = \int_0^{r} \frac{1}{V(s)}\;.
\end{equation}
Here $V$ is parameterized with respect to the distance $s$ at time $t$ of the wave front from its origin.
As shown in the proof of \cite[Theorem 3.1]{KowSchBon10}, the
fact that $\beta^*$ is of standard form together with~(\ref{eq:lem1Tc0})
implies that there exist a constant $c$ such that
$T(r) = r/c$ for all $r >0$.
But then from ~(\ref{eq:lem1Tc0}) it follows $0<r/c_0\leq r/c<\infty$ and consequently
$0 < c \leq c_0 < \infty$.
\end{proof}

Finally we explain the above notation for the standard wave
equation:
\begin{remark}
In the case of the standard wave equation the wave front is the
support of the Green function $\green_0$, the wave front velocity is
$c_0$, and $T(\abs{\vx})=\frac{\abs{\vx}}{c_0}$ denotes the travel
time of the wave front.
\end{remark}

\section{Strong Causality of Attenuation Laws}
\label{sec:n}
In this section we analyze causality properties of attenuation laws. We split the section
into two parts, where the first concerns numerical studies to determine the kernel function
$K$, defined in \req{kernel}, and the second part contains analytical investigations.

In Figures~\ref{fig:powlaw},~\ref{fig:szabo} and~\ref{fig:thviscous} we represent the
attenuation kernels according to power, Szabo's, and the thermo-viscous law.
\begin{figure}[htb]
\begin{center}
\includegraphics[height=4.0cm,angle=0]{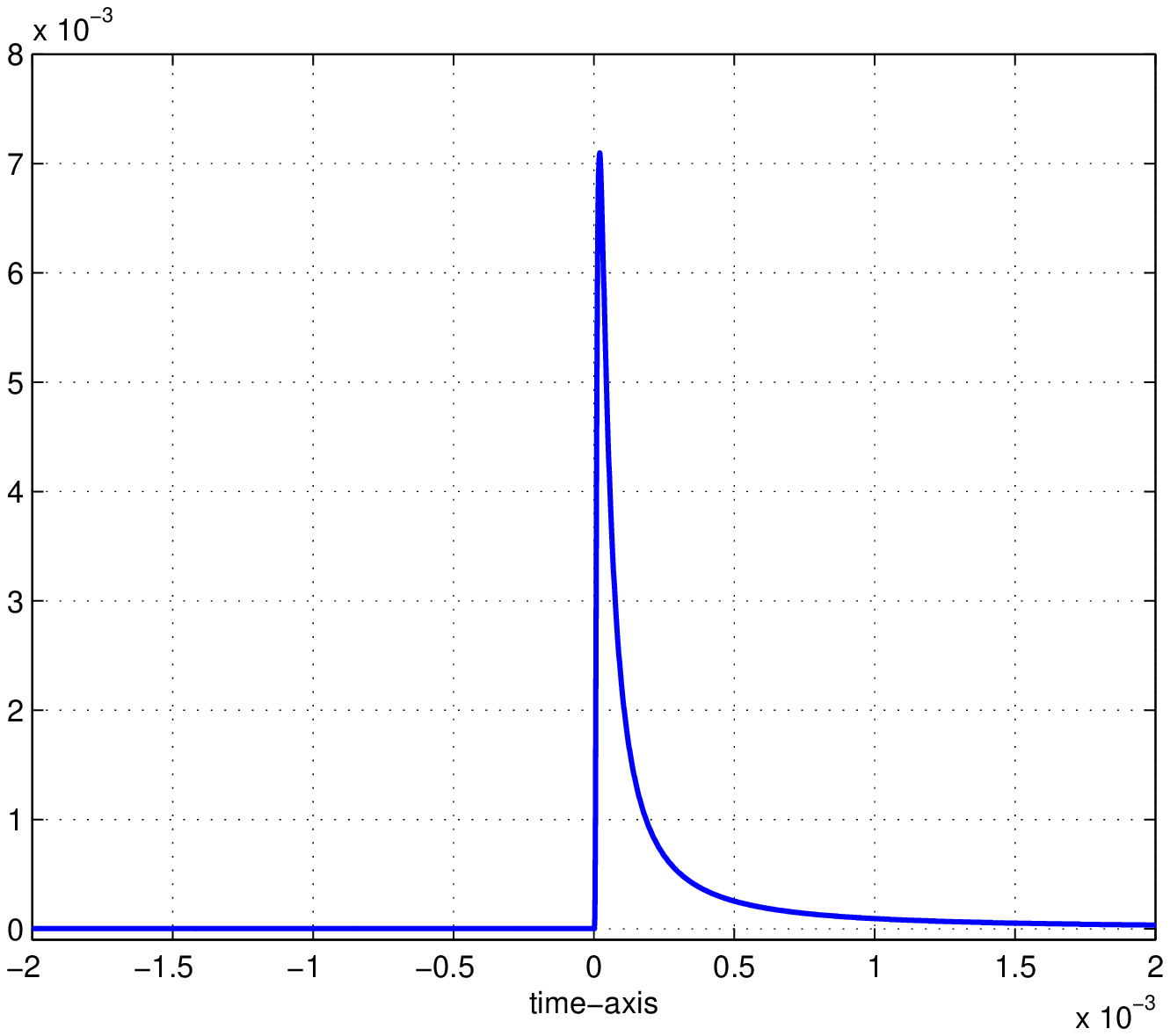}
\includegraphics[height=4.0cm,angle=0]{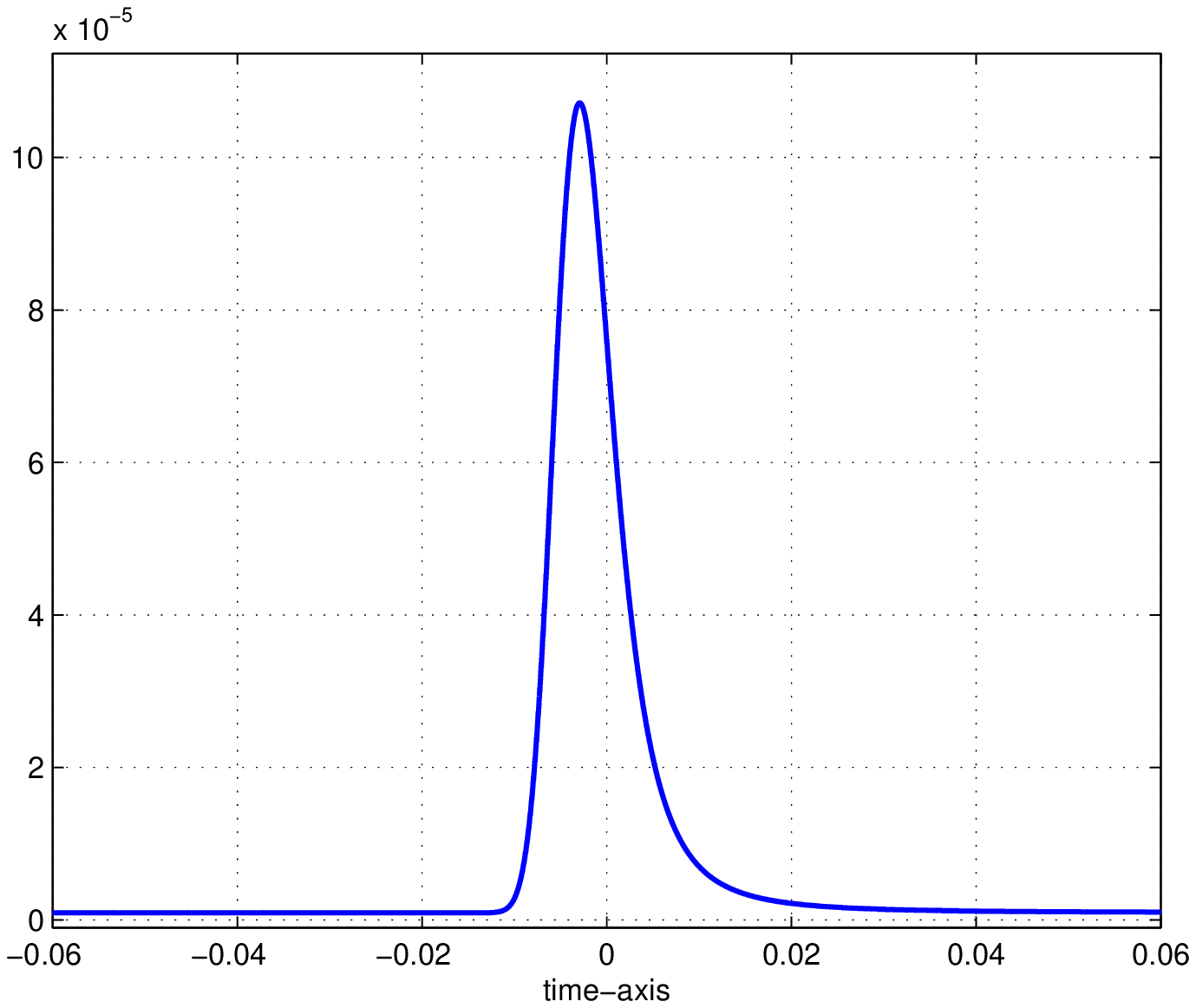}\\
\includegraphics[height=4.0cm,angle=0]{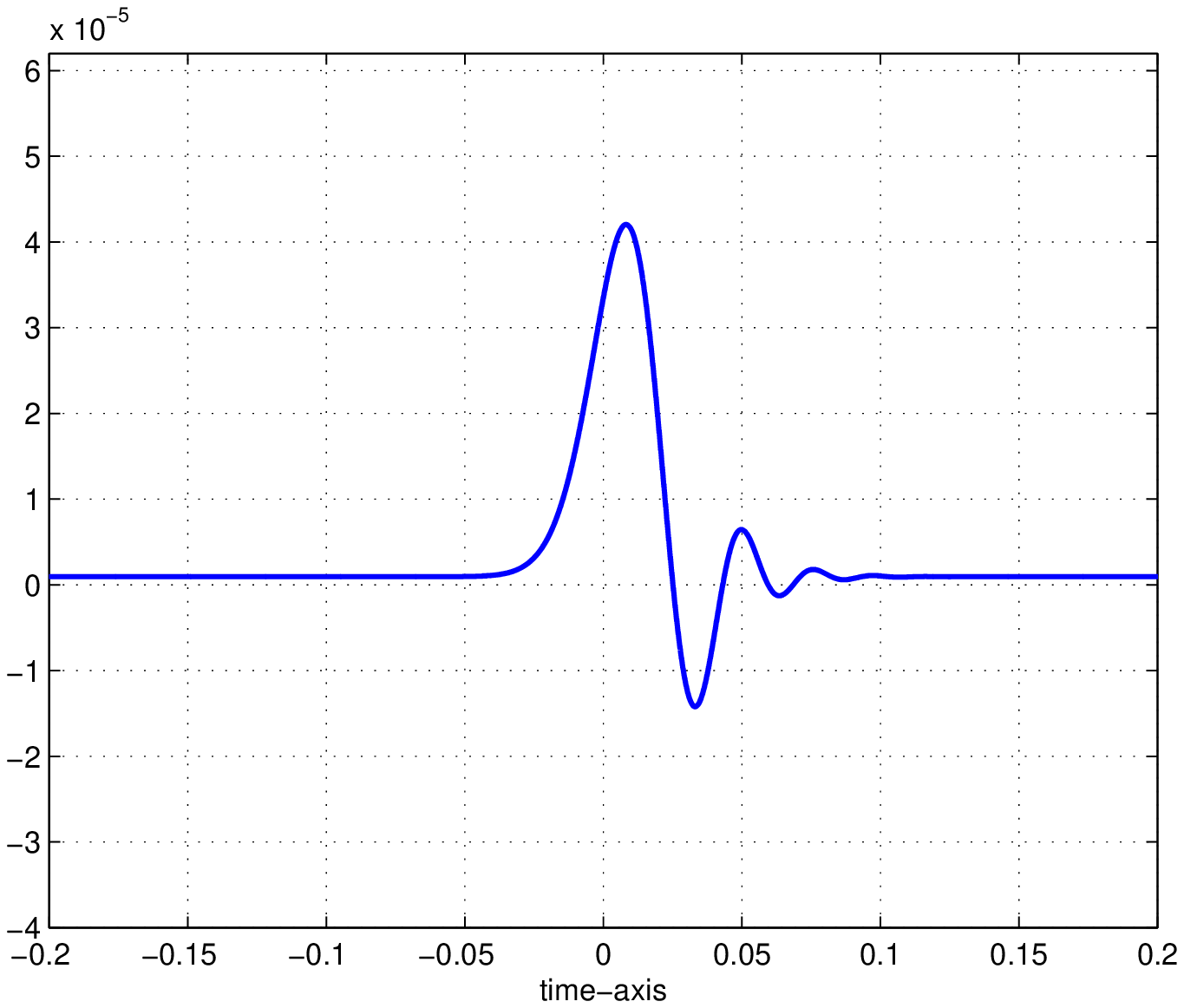}
\includegraphics[height=4.0cm,angle=0]{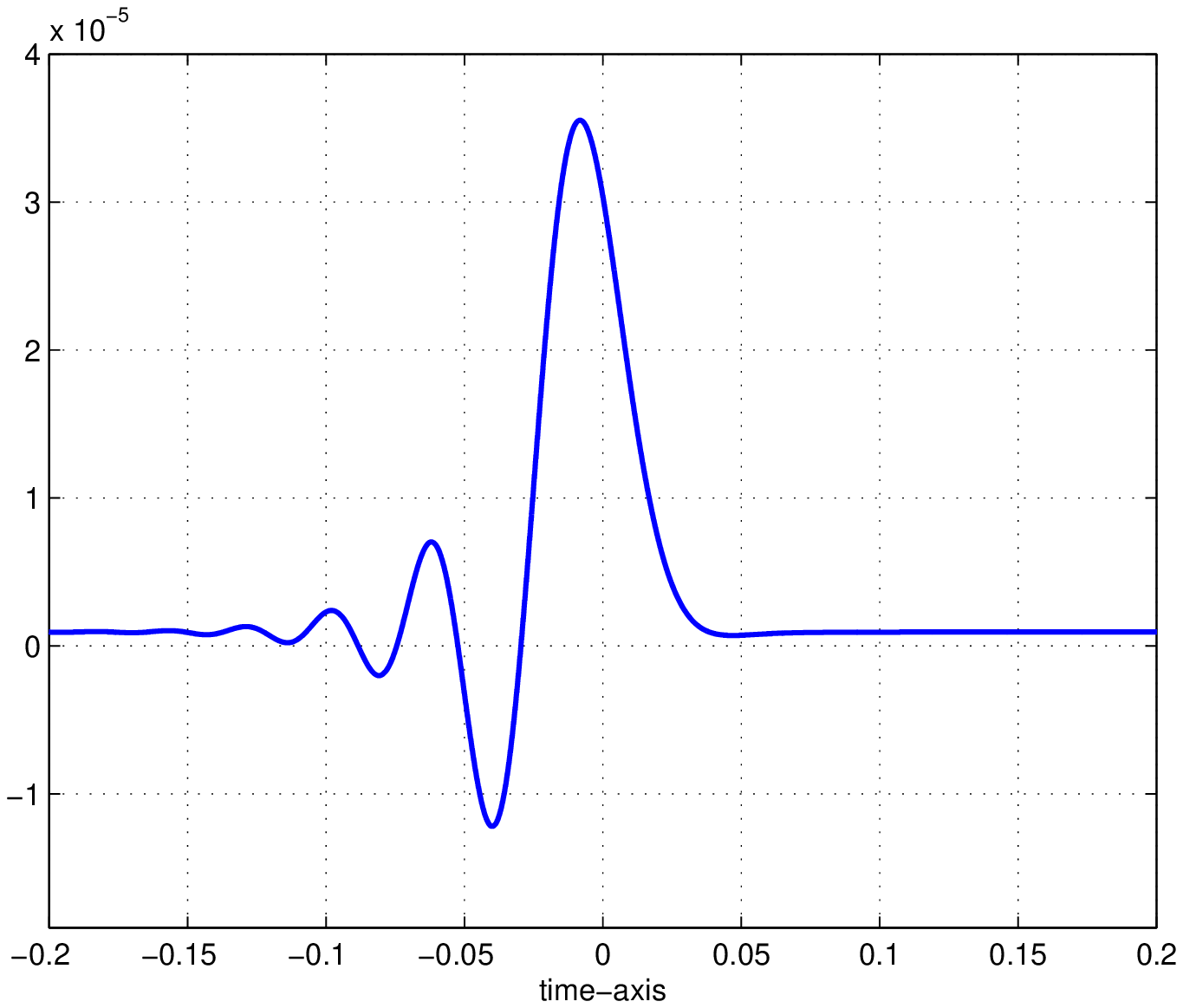}
\end{center}
\caption{Simulation of $K(\vx,t)$ for the frequency power law with $(\gamma,\alpha_0)\in
\set{(0.5,0.1581),\,(1.5,0.0316),\,(2.7,0.0071),\,(3.3,0.0027)}$, $c_0=1$ and
$\abs{\vx}=\frac{1}{4}$. In the first example $\gamma < 1$ and thus the function is causal.
For all other cases it is non causal.} \label{fig:powlaw}
\end{figure}
\begin{figure}[htb]
\begin{center}
\includegraphics[height=4.0cm,angle=0]{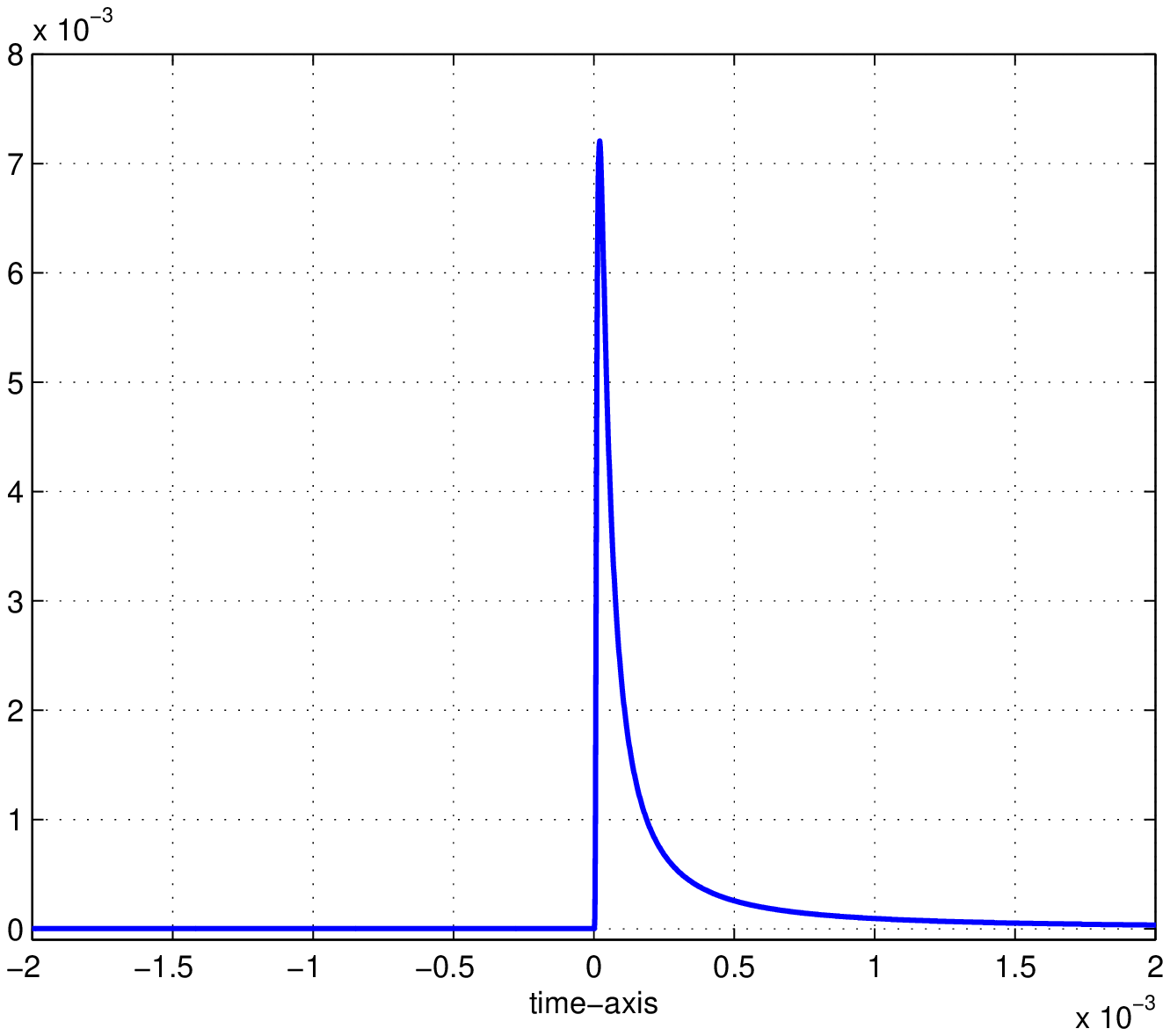}
\includegraphics[height=4.0cm,angle=0]{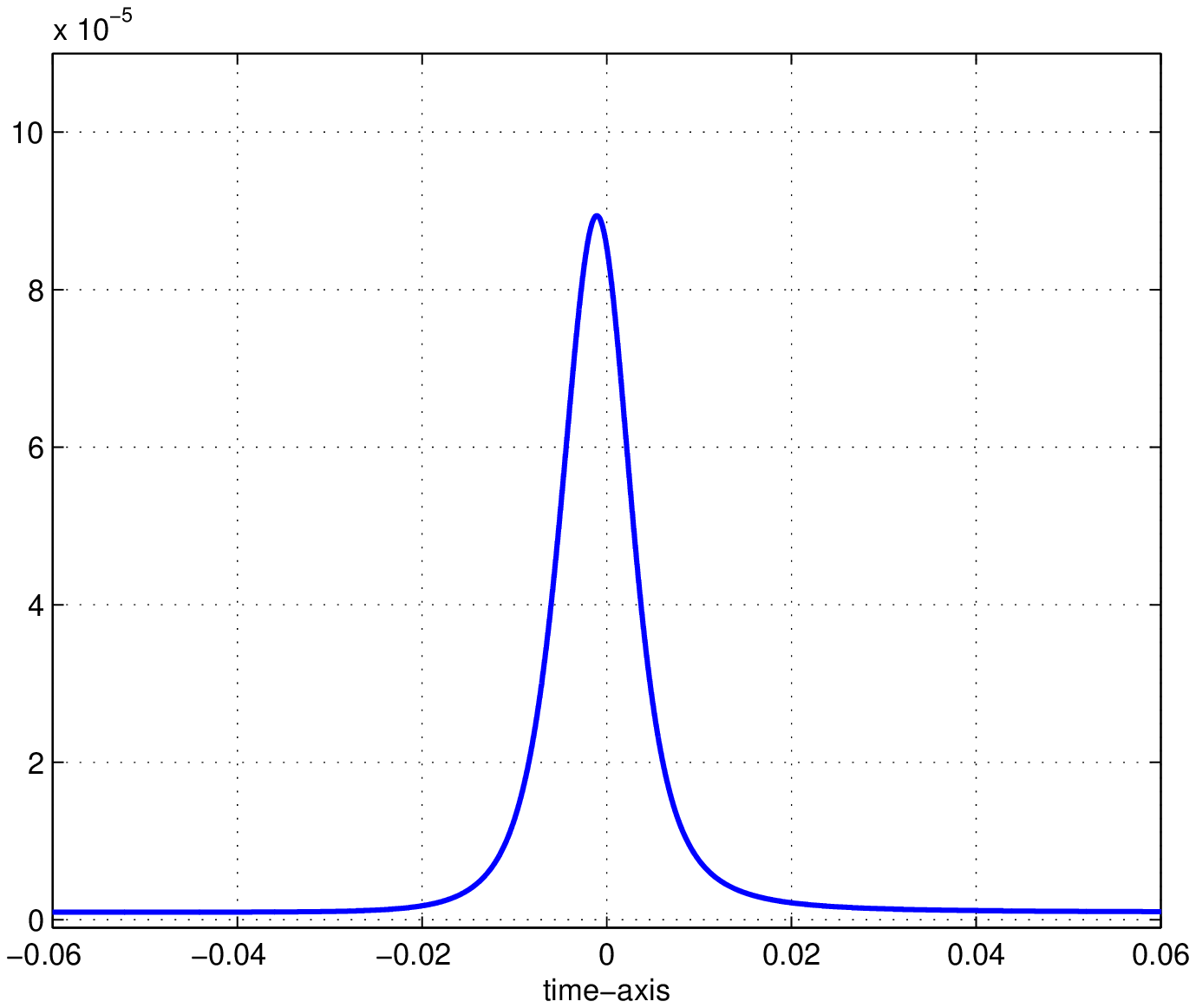}\\
\includegraphics[height=4.0cm,angle=0]{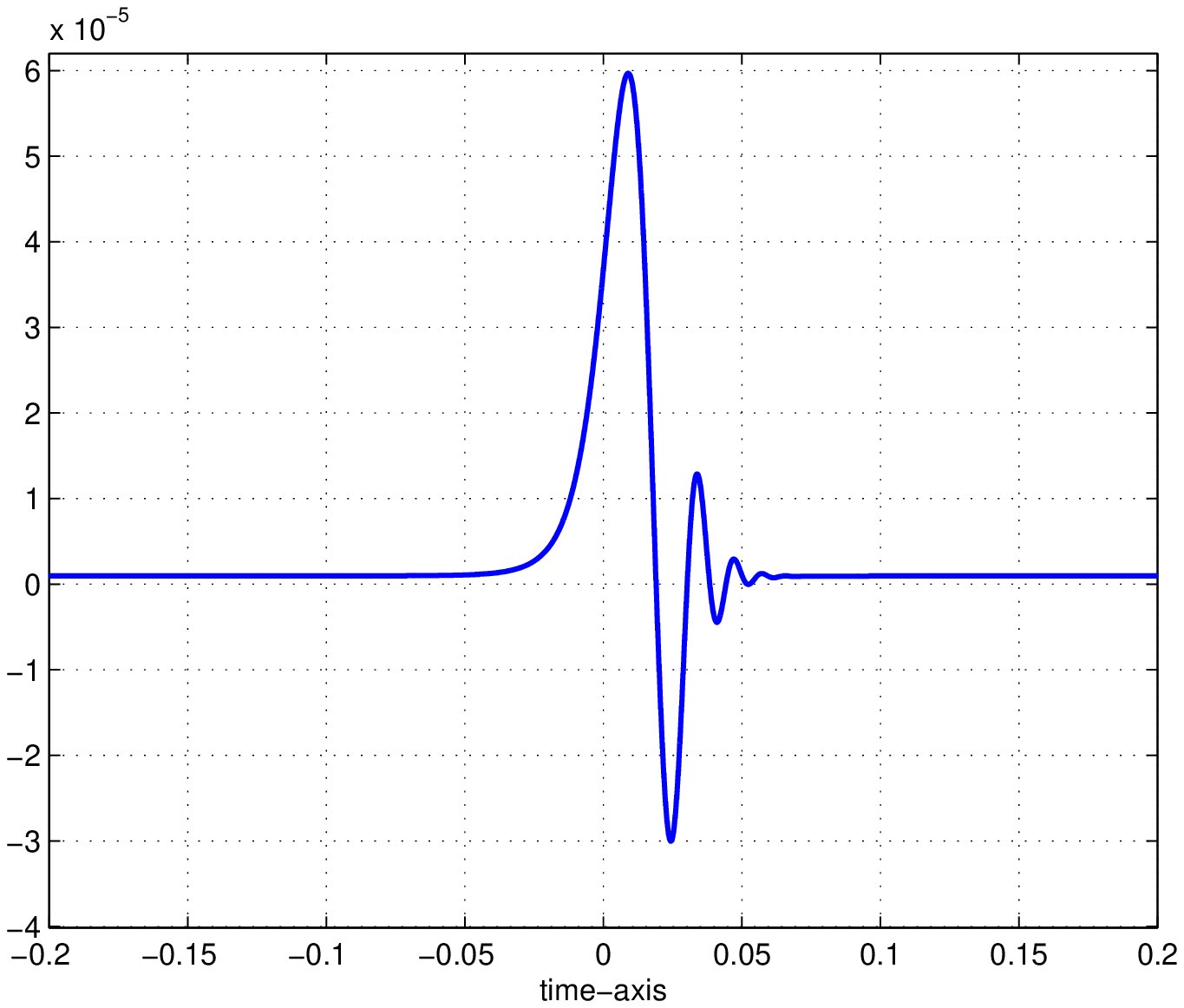}
\includegraphics[height=4.0cm,angle=0]{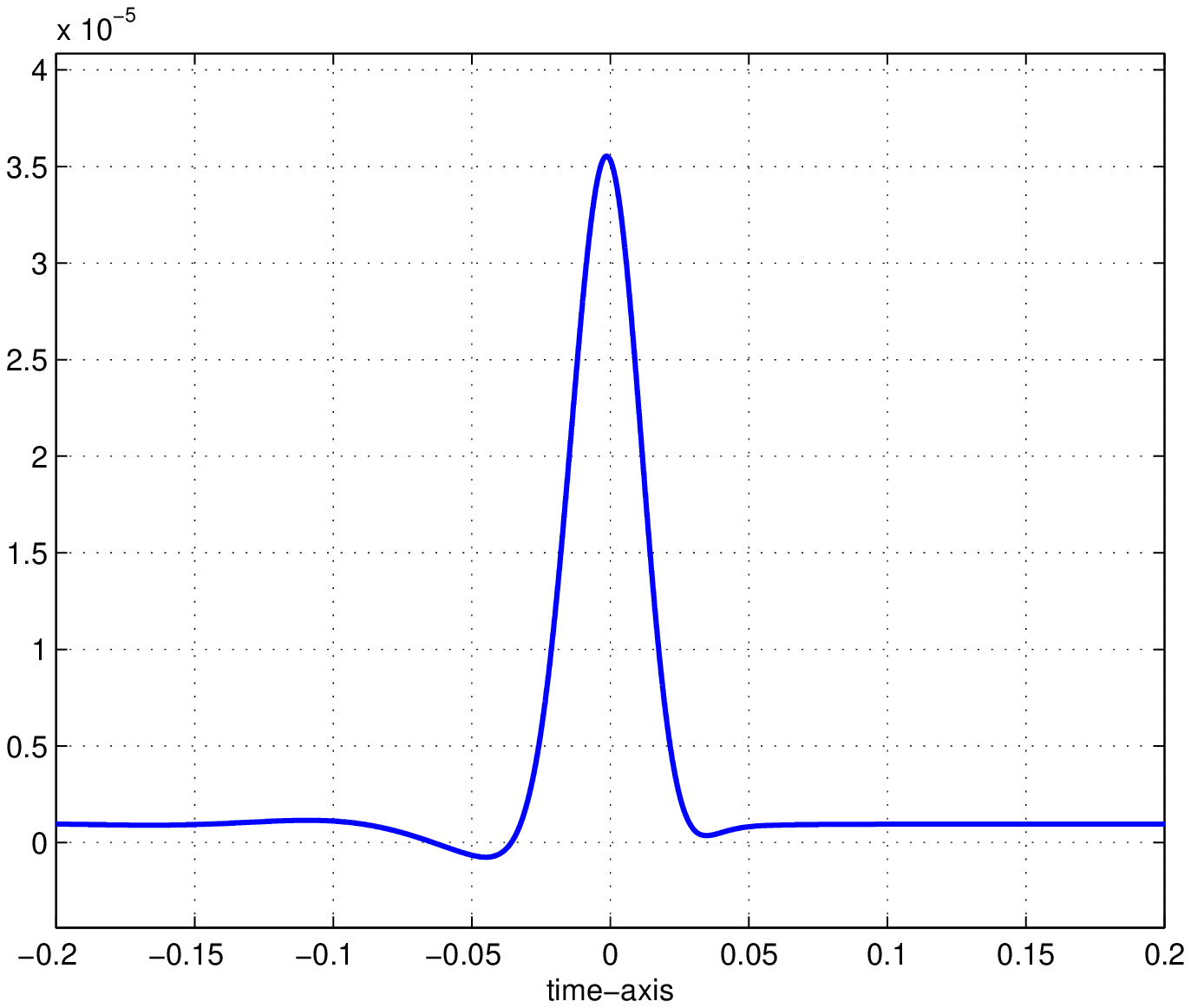}
\end{center}
\caption{Simulation of $K(\vx,t)$ for Szabo's frequency law with $(\gamma,\alpha_0)\in
\set{(0.5,0.1581),\,(1.5,0.0316),\,(2.7,0.0071),\,(3.3,0.0027)}$, $c_0=1$ and
$\abs{\vx}=\frac{1}{4}$. } \label{fig:szabo}
\end{figure}
\begin{figure}[htb]
\begin{center}
\includegraphics[height=4.0cm,angle=0]{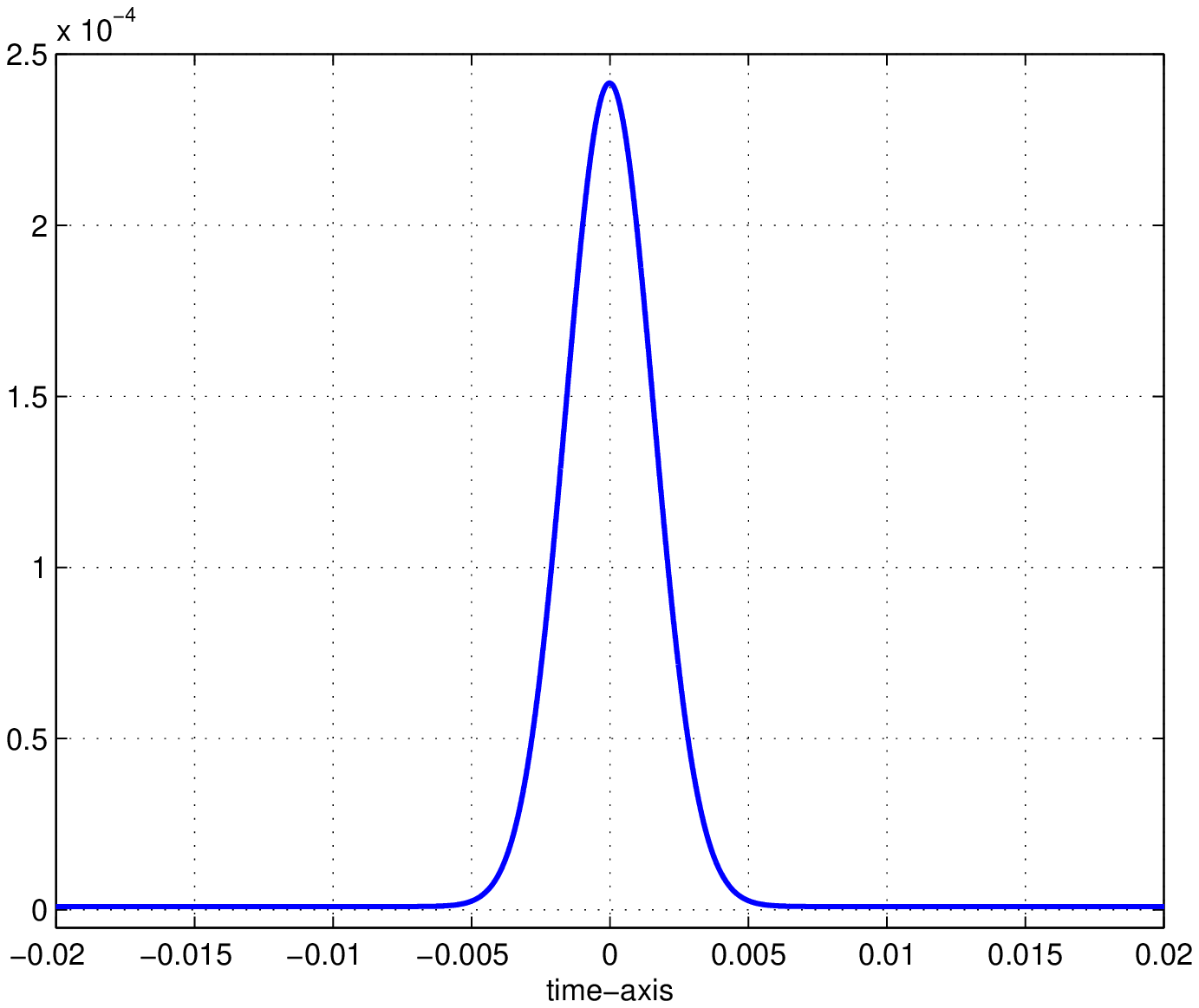}
\includegraphics[height=4.0cm,angle=0]{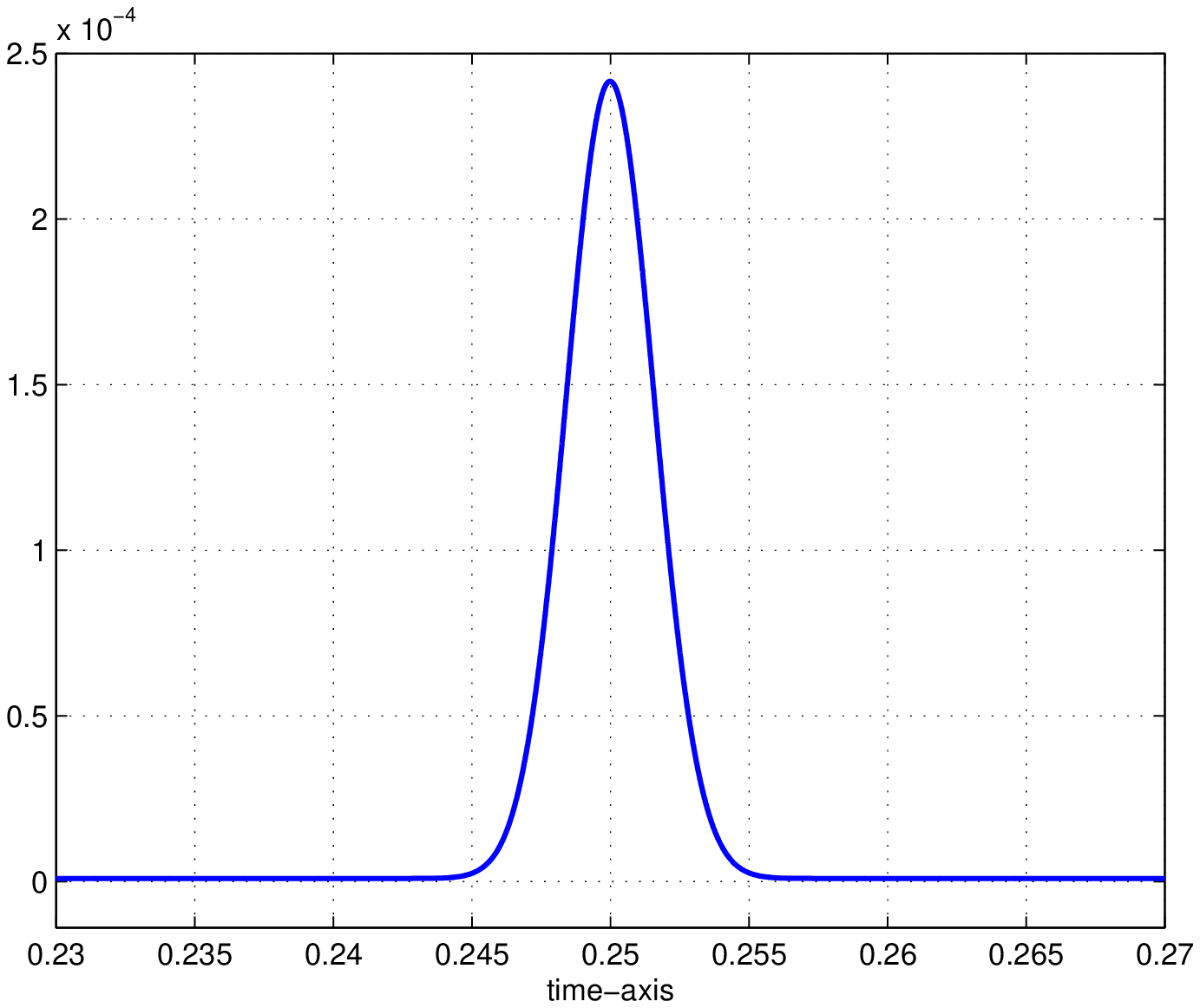}
\end{center}
\caption{{\sl Left:} $K(\vx,t)$ defined by the complex thermo-viscous attenuation law with
$\tau_0=10^{-5}$, $c_0=1$ and fixed $\abs{\vx}=\frac{1}{4}$. {\sl Right:} The proposed law~(\ref{eq:powlaw2})
for $\gamma=2$ with $\alpha_1=1$, $\tau_0=10^{-5}$, $c_0=1$ and fixed $\abs{\vx}=\frac{1}{4}$ is causal. } \label{fig:thviscous}
\end{figure}
The figures already indicate that power laws with index greater than $1$ violate causality.
In the following we support these computational studies by analytical considerations.
Thereby we make use of distribution theory, which we recall first. Generally speaking
\emph{Distributions} are generalized functions:
\begin{definition}
\label{def:distribution}
We use the abbreviations:
\begin{itemize}
\item $\D:=C_0^\infty(\R,\C)$ is the space of infinitely often differentiable functions
from $\R$ to $\C$ which have compact support. \index{$\D$}
\item $\S$ is the space of infinitely often differentiable functions from $\R$ to $\C$
      which are \emph{rapidly decreasing}. \index{$\S$}
      A function $f: \R \to \C$ is rapidly decreasing if for all $i,j \in \N_0$
      \begin{equation*}
      \abs{x}^i \abs{f^{(j)}(x)} \to 0 \text{ for } \abs{x} \to \infty\;.
      \end{equation*}
      \index{Function, rapidly decreasing}
\item $\S$ is a locally convex space (see \cite{Yos95} for a definition) with the topology
      induced by the family of semi-norms
      \begin{equation*}
      p_{P,j}(f) = \sup_{x \in \R} \abs{P(x) f^{(j)}(x)} ,
      \end{equation*}
      where $P$ is a polynomial and $j \in \N_0$. The topology on a locally convex set is defined as follows:
      $U \subseteq \S$ is open, if for every $f \in U$ there exists $\ve > 0$ and a finite non-empty set $J'$
      of polynomials and a finite set of indices $K'$ such that
      \begin{equation*}
      \bigcap_{P \in J', k \in K'} \set{g \in \S: p_{P,k}(g-f) < \ve} \subseteq U\;.
      \end{equation*}
\item The space of \emph{tempered distributions}, $\S'$, is the space of linear continuous functionals
      and $\S$. \index{tempered distributions, $\S'$}
\item A functional $L: \S \to \C$ is continuous if there exists a constant $C > 0$ and a seminorm $p_{P,j}$ such that
      \begin{equation*}
      \abs{L u} \leq C p_{P,j}(u)\,, \text{ for every } u \in \S
      \end{equation*}
      (see \cite[Sect. I.6, Thm. 1]{Yos95}
\end{itemize}
\end{definition}

In the following we give some examples of tempered distributions and review some of their properties.
The examples are taken from \cite[Sec6.2, Ex. 3]{Yos95} and \cite[Remark 6]{DauLio02_2}.
\begin{example}{\bf (Examples of Tempered Distributions)}
\label{ex:distr_ex}
\begin{itemize}
\item Let $1 \leq p \leq \infty$ and $f \in L^p(\R,\C)$, then the linear operator $T\phi = \int_{\R} f(x)\phi(x)\,dx$ is
a tempered distribution. In the following we identify $f$ and $T$, and this clarifies the terminology $f \in \S'$ later on.
\item $\S \subseteq \S'$ - thereby already the above relation between functions and tempered distributions is used.
\item Distributions with compact support are tempered distributions. For instance the $\delta$-Distribution is a tempered distribution.
\item Polynomials are tempered distributions.
\item The functions $f$ of $L_{loc}^1(\R)$
which are uniformly bounded by a polynomial for $\abs{x}$ sufficiently large, are tempered distributions.
\footnote{A function $f$ is an element of $L_{loc}^1(\R)$ if it is in $L^1$ on every compact set.}
\end{itemize}
\end{example}

\begin{lemma}
\label{le:distr_ex}
\begin{itemize}
\item The pointwise limit $f:\R \to \C$ of a sequence of functions $\set{f_n:\R \to \C} \subseteq \S'$,
      is again a tempered distribution.
\item Let $f \in \S'$, then $\fourier{f} \in \S'$ and $\ifourier{f} \in \S'$.
\end{itemize}
\end{lemma}

In the following we review Theorem 4 on p294 ff from \cite{DauLio02_5} which characterized
when a generalized function $f\in\S'(\R)$ is causal, that is, when
$\supp(f) \subseteq [0,\infty)$. Below we use the following notation
\begin{equation*}
\C_\ve:=\set{z\in\C : \Im(z)\geq \ve}\;.
\end{equation*}
\index{$\C_\ve$}
\begin{theorem} (Theorem 4 on p294 ff in \cite{DauLio02_5})
\label{th:lion}
Let $f \in \S'(\R)$. Then $f$ is causal
\footnote{In Theorem 4 on p294 ff in \cite{DauLio02_5} the assumption that $f$ is strongly causal is
expressed by $f \in \mathring{\D}_+$, which is the set of distributions with support in $[0,+\infty)$.}
if and only if
\begin{enumerate}
\item \label{it1_Lions} There exists a function $F:\C_0:=\set{\xi + \i \eta : \eta \geq 0}\to\C$, which is holomorphic in the
interior $\mathring{\C}_0:=\set{\xi + \i \eta : \eta > 0}$.
\index{$\C_0$}
\footnote{A function $F:\C_0 \to\C$ is holomorphic in $\mathring{\C}_0$ if it is complex
differentiable in $\mathring{\C}_0$. Sometimes the functions are also
refered to as analytic or regular functions or conformal maps.}
\item \label{it2_Lions}
For all fixed $\eta>0$ and $\xi\in\R$, $F(\xi + \i \eta)$ is a tempered distribution
with respect to the variable $\xi$ and for $\eta\to 0$ $F(\xi + \i \eta)$ is convergent
(with respect to the weak topology on $\S'$) to $\fourier{f}:\R\to\C$.
\footnote{A function $F:\C_0\to\C$ which satisfies Items \ref{it1_Lions},\ref{it2_Lions} of Theorem \ref{th:lion} is called
\emph{holomorphic extension} of $\fourier{f}$.}
\item \label{it3_Lions} For every $\ve > 0$, there exists a polynomial $P$ such that
\begin{equation*}
\abs{F(z)} \leq P(\abs{z})  \ltext{for} z \in \C_\ve\,.
\end{equation*}
\end{enumerate}
\end{theorem}

\begin{remark}
\label{re:holomophic}
The definition of the Fourier transform in this chapter has a different sign as
in~\cite{DauLio02_5} and consequently also $\C_0$ denotes the upper half plane and not the
lower half plane as in \cite{DauLio02_5}.
\end{remark}

For analyzing attenuation laws, we use the following corollary, which is derived from Theorem
\ref{th:lion}.
\begin{corollary}
\label{co:lion}
Let $\alpha^*: \R \to \C$ be continuous and let there exist a holomorphic extension to $\C_0$, which for the sake of
simplicity of notation is again denoted be $\alpha^*$. In addition, let $\alpha^* (\xi+\i \eta) \to \alpha^* (\xi)$ for $\eta \to 0$
pointwise. We denote by $\alpha: \C_0 \to \C$ the real part
of $\alpha^*$. \footnote{$\alpha: \C_0 \to \C$ extends the function $\omega \in \R \to \alpha(\omega)$ but is not an
holomorphic extension.}
Moreover, we assume that there exists a constant $C$ such that
\begin{equation}
\label{estimate:rapidly}
\alpha(\omega) \geq C \text{ for all }\omega \in \R\;.
\end{equation}

\begin{enumerate}
\item If in addition
\begin{equation}
\label{estimate:polynomial}
\alpha(z) \geq C \text{ for all } z \in \C_0\;.
\end{equation}
Then, for every $\vx \in \R^3$, the function
\begin{equation}
\label{31}
t \to K(\vx,t) := \frac{1}{\sqrt{2\pi}}\ifourier{\ef{-\alpha^*(\cdot)\abs{\vx}}}(t)
\end{equation}
is causal.
\item On the other hand, if there exists $C_1>0$, $\mu>0$ and $C_2 \in \R$ and a sequence $\set{z_n}$ in $\mathring{\C}_0$ such that
\begin{equation}
\label{estimate:polynomial_negativ}
\alpha(z_n) \leq - C_1 \abs{z_n}^\mu - C_2\,,
\end{equation}
then $K$ violates causality.
\end{enumerate}
\end{corollary}

\begin{proof}
Let $\vx \in \R^3$ fixed. We apply Theorem \ref{th:lion} to $f(\cdot)=K(\vx,\cdot)$. Therefore, we have
\begin{equation*}
\sqrt{2\pi} \fourier{f}(\omega) = \ef{ -\alpha^*(\omega) \abs{\vx}}\;.
\end{equation*}
Under the assumption (\ref{estimate:rapidly}), taking into account that $\alpha_{pl}^*$ is continuous, $\fourier{f}$ is in $L_{loc}^1(\R)$
and bounded by a constant polynomial, thus in $\S'$ (cf. Example \ref{le:distr_ex}), and consequently, according to Lemma \ref{le:distr_ex},
$f \in \S'$. Therefore, the general assumption of Theorem \ref{th:lion} is satisfied.

The function $z \in \C_0 \mapsto F(z):=\ef{ -\alpha^*(z) \abs{\vx}}$ is an extension of
$2\pi \fourier{f}(\omega)$, which is holomorphic in $\mathring{\C}_0$.
Thus Item \ref{it1_Lions} of Theorem \ref{th:lion} holds.

\begin{itemize}
\item For proving the first assertion, it follows from (\ref{estimate:polynomial}) that for
$z = \xi + \i \eta$ with $\eta \geq 0$,
\begin{equation}
\label{eq:new}
\abs{F(z)} \leq \ef{-C\abs{\vx}}\,,
\end{equation}
which, in particular, shows that for all $\eta \geq 0$, the functions $\xi \to F(\xi+\i\eta)$ is a tempered
distribution (cf. Lemma \ref{le:distr_ex}).
Hence Item~\ref{it1_Lions} of Theorem \ref{th:lion} holds.

Moreover, since by assumption $\alpha^*(\xi+\i \eta) \to \alpha^*(\xi)$ for $\eta \to 0$ pointwise, $F(\xi+\i\eta)$ converges to
$F(\xi)=2\pi \fourier{f}(\xi)$ pointwise.
Because the limit is a tempered distribution and the convergence is with respect to the weak topology $\S'$
(cf. Lemma \ref{le:distr_ex}).
Hence Item~\ref{it2_Lions} of Theorem \ref{th:lion} holds.

Moreover, from (\ref{eq:new})
it follows that $\abs{F(z)}$ is bounded by a constant polynomial.
Hence Item~(\ref{it3_Lions}) of Theorem~\ref{th:lion} holds and
therefore Theorem \ref{th:lion} guarantees that $t \mapsto K(\vx,t)$ is causal.

\item
For the second case, Item \ref{it3_Lions} of Theorem \ref{th:lion} is violated. Consequently, $K(\vx,\cdot)$ is not causal.
\end{itemize}
\end{proof}

\subsection*{Power Laws}

\begin{theorem}\label{coro:powlaw1}
Let $0 < \gamma \in \R$, be not an odd number, and $\omega \in \R \mapsto \alpha_{pl}^*(\omega)=\ta_0 (-\i\omega)^\gamma$ be
the power law attenuation coefficient from (\ref{eq:powlaw1}), with
$\ta_0=\alpha_0/\cos \left(\frac{\pi}{2}\gamma\right)$ as in (\ref{ta0}). Then, the function $K$, defined in
(\ref{31}), is causal if and only if $\gamma\in (0,1)$.
\end{theorem}

\begin{proof}
Let $\vx \in \R^3$ be fixed. The function $z \in \C \mapsto \alpha_{pl}^*(z)=\ta_0 (-\i z)^\gamma$ is the
holomorphic extension of $\omega \in \R \mapsto \alpha_{pl}^*(\omega)$. We prove or disprove causality
by using Corollary \ref{co:lion}.

For $z=\abs{z} \ef{\i\phi}$ it follows from (\ref{eq:powlaw1}) that
$$
    \Re((-\i z)^\gamma)
        = \Re(\abs{z}^\gamma\,\ef{\i \gamma\,(\phi-\pi/2)})
        = \abs{z}^\gamma\,\cos(\gamma(\phi-\pi/2))\,.
$$
This implies that
\begin{equation}
\alpha_{pl}(z) =  \ta_0 \Re((-\i z)^\gamma)  =  \ta_0 \cos(\gamma(\phi-\pi/2))\abs{z}^\gamma\;.
\end{equation}
In particular, if $z = \omega \in \R$, then $\phi$ is either $0$ or $\pi$. Taking into account the definition of $\ta_0$
and that the $\cos$-function is symmetric around the origin, it follows
that
\begin{equation}\label{relalphapowlaw1a}
\alpha_{pl}(\omega) = \alpha_0 \abs{w}^\gamma \geq 0\;.
\end{equation}
Thus (\ref{estimate:rapidly}) holds.

\begin{itemize}
\item Let $\gamma \in (0,1)$: Every $z=\abs{z}\,\ef{\i\phi}\in \C_0$
satisfies $\phi\in [0,\pi]$. Consequently $\gamma(\phi-\pi/2) \in [-\pi/2,\pi/2]$ and thus
$\cos(\gamma(\phi-\pi/2))$ is uniformly non-negative.
Even more for $\gamma \in (0,1)$ the coefficient $\ta_0$, defined
in (\ref{ta0}), is positive. In summary, we have that there exists a constant $C_1 \geq 0$ such that
\begin{equation}\label{relalphapowlaw1b}
\alpha_{pl}(z) \geq C_1 \abs{z}^\gamma \geq 0 \ltext{ for } z\in \C_0\;.
\end{equation}
Thus (\ref{estimate:polynomial}) holds and application of Corollary \ref{co:lion} shows that $K$ is causal.

\item Let $\gamma\in (1,3)\cup (5,7)\cup\cdots$. Then $\ta_0<0$. The sequence
\begin{equation}
\label{eq:34a}
\set{z_n := n \ef{\i \pi/2} = \i n}_{n \in \N}
\end{equation}
consists of elements of $\mathring{\C}_0$ and
satisfies assumption (\ref{estimate:polynomial_negativ}), that is,
\begin{equation}
\label{eq:azn1}
\alpha_{pl}(z_n) = \underbrace{\ta_0}_{<0} \abs{z_n}^\gamma\;.
\end{equation}
Application of Corollary \ref{co:lion} shows that $K$ is not causal.

\item  Let $\gamma\in (3,5)\cup (7,9)\cup\cdots$. We fix some $0<\delta<\pi/2$, and define
\begin{equation*}
\phi := \left(1+\frac{1}{\gamma} \right)\frac{\pi}{2}  + \frac{\delta}{\gamma}\;.
\end{equation*}
The sequence
\begin{equation}
\label{eq:34b}
\set{z_n:=n\,\ef{\i \phi}}
\end{equation}
consists of elements of $\mathring{\C}_0$. Under the above assumptions, it follows that $\ta_0 > 0$
and therefore
\begin{equation}
\label{eq:azn2}
\alpha_{pl} (z_n) = \underbrace{\ta_0}_{>0} \underbrace{\cos(\pi/2+\delta)}_{<0} \abs{z_n}^\gamma\;.
\end{equation}
Thus form Corollary \ref{co:lion} the assertion follows.
\end{itemize}
\end{proof}

In the following we analyze the following family of variants of power laws:
\begin{equation}\label{szabo_model}
\alpha_{pl+}^*(\omega) = \ta_0 (-\i \omega)^\gamma  + \alpha_1(-\i \omega)\,,
\end{equation}
which have been considered in \cite{Sza95,WatHugBraMil00}.

\begin{theorem}
\label{le:pow_law_extended}
Let $0 < \gamma \notin \N$ and $\alpha_{pl+}^*$ as defined in (\ref{szabo_model}). Moreover, let $K$
be as in (\ref{31}).
Then, if $\gamma >1$, $K$ is not causal.
For $\gamma \in (0,1)$ $K$ is causal if and only if $\alpha_1\in [0,\infty)$.
\end{theorem}

\begin{proof}
The holomorphic extension of $\omega \in \R \to \alpha_{pl+}^*(\omega)$ is the function
\begin{equation*}
\alpha_{pl+}^*(z) = \ta_0 (-\i z)^\gamma  + \alpha_1(-\i z)
\end{equation*}
and consequently
\begin{equation*}
\alpha_{pl+}(z)
= \ta_0 \abs{z}^\gamma \left(1+\frac{\alpha_1}{\ta_0} \abs{z}^{1-\gamma} \right)
\cos \left(\gamma \left(\phi-\frac{\pi}{2}\right)\right)\;.
\end{equation*}

\begin{itemize}
\item For $\gamma > 1$ we have that $1+\frac{\alpha_1}{\ta_0} \abs{z}^{1-\gamma} \to 1$ for
$\abs{z} \to \infty$.
Let $\set{z_n}$ as defined in (\ref{eq:34a}) or (\ref{eq:34b}). Then, since for both
sequences $\abs{z_n} \to \infty$, it follows from (\ref{eq:azn1}), (\ref{eq:azn2}) that the according
sequences $\set{z_n}$ satisfy (\ref{estimate:polynomial_negativ}) for $n$ sufficiently large, respectively.
Thus $K$ is not causal.

\item For $\gamma \in (0,1)$ and $\alpha_1 \geq 0$ the assertion follows already from the fact
that $\alpha_{pl+}(\omega) \geq \alpha_{pl}(\omega)$ and that the later already satisfies (\ref{estimate:polynomial}).
Thus $K$ is causal.

\item Let $\gamma \in (0,1)$ and $\alpha_1 < 0$. Then for some $0 < \delta < -\alpha_1$ fixed, we can find
a constant $C_2$ such that for all $z \in \C$
\begin{equation*}
\ta_0 \abs{z}^\gamma + \alpha_1 \abs{z} \leq \underbrace{(\alpha_1+\delta)}_{< 0} \abs{z} - C_2\;.
\end{equation*}
Consequently, for $\set{z_n=\i n}$, we have
\begin{equation*}
\alpha_{pl+}(z_n) \leq (\alpha_1+\delta) \abs{n} - C_2\;.
\end{equation*}
which shows (\ref{estimate:polynomial_negativ}). Thus $K$ is not causal.
\end{itemize}
\end{proof}

\subsection*{Powerlaw with $\gamma = 1$}

\begin{theorem} \label{thcaus02}
Let $\alpha_{pl}^*$ be as in defined in (\ref{alphagamma1}).
Then the function $K$, defined in (\ref{31}) is not causal.
\end{theorem}

\begin{proof}
First, we prove that
\begin{equation*}
\begin{aligned}
z \in \C_0 \mapsto \hat{\alpha}_{pl}^*(z) := \alpha_0 z
        + \i \frac{2 \alpha_0}{\pi} z \log \left(\frac{z}{\omega_0}\right)
\end{aligned}
\end{equation*}
is the holomorphic extension of $\omega \to \alpha_{pl}^*(\omega)$.
This assertion follows from the facts
\begin{equation*}
\begin{aligned}
\hat{\alpha}_{pl}^*(\omega)& =\alpha_{pl}^*(\omega)  \text{ for } & \omega>0\,,\\
\lim_{\eta \to 0+} \i \frac{2}{\pi} \log \left(\frac{\omega + \i \eta}{\omega_0}\right)
      & = \i \frac{2}{\pi}\,\log\left|\frac{\omega}{\omega_0}\right|
        - 2 \text{ for } & \omega<0\;.
\end{aligned}
\end{equation*}
Since
$$
   \alpha_{pl}(\omega) = \Re( \alpha_{pl}^*(\omega)) = \alpha_0\,\abs{\omega} \geq 0\,,
$$
Corollary \ref{co:lion} is applicable.
For the elements of the sequence $\set{z_n:=\i n}_{n \in \N}$ in $\C_0$
$$
  \hat{\alpha}_{pl}^*(z_n)
       = - \frac{2 \alpha_0}{\pi} n \log \left(\frac{n}{\omega_0}\right)
$$
is real and therefore equals $\alpha_{pl}(z_n)$ and thus (\ref{estimate:polynomial_negativ}) holds.
Thus Corollary \ref{co:lion} gives the assertion.
\end{proof}

\subsection*{Szabo's Model:}

\begin{proposition}\label{theo:szabo}
For $\alpha_0>0$ let $\alpha_{sz}^*$ be the coefficient of Szabo's model~(\ref{alpha*szabo}).
Then, for $\gamma\in (0,1)$, the function $K$ (\ref{eq:kernel}) is a causal function and for
$\gamma>1$ with $\gamma\not\in\N$, $K$ violates causality.
\end{proposition}

\begin{proof}
Without loss of generality we assume that $c_0=1$.
The holomorphic extension of $\alpha_{sz}^*:\R \to \C$ from (\ref{alpha*szabo}) is
\begin{equation*}
 z\in \C_0 \to \alpha_{sz}^*(z) = (-\i z)
  \left[ \sqrt{1 + 2 \ta_0 (-\i z)^{\gamma-1}} - 1\right]\,.
\end{equation*}

First, we make some general manipulations which can be used in several ways:
Let $z = \xi+\i \eta \in \C_0$. We use the polar representation
\begin{equation*}
z = \abs{z} \ef{\i \phi}\,, \quad \phi := \phi(z) \in \left[0,\pi\right]\;.
\end{equation*}
Then
\begin{equation*}
\alpha_{sz}^*(z) = \abs{z} \ef{\i (\phi-\pi/2)} \Psi(z)\,,
\end{equation*}
where
\begin{equation}
\label{eq:delta}
\begin{aligned}
\Psi(z) &:= \sqrt{\hat{\Psi}(z)} - 1\,,\\
\hat{\Psi}(z) &:= 1 + 2 \ta_0 \abs{z}^{\gamma-1} \ef{\i \delta}\,,\\
\delta &:= \delta(z):=(\phi-\pi/2)(\gamma-1)\;.
\end{aligned}
\end{equation}
With this notation we have
\begin{equation}
\label{eq:hatpsi}
\begin{aligned}
\Re (\hat{\Psi}(z)) &= 1 + 2 \ta_0 \cos(\delta) \abs{z}^{\gamma-1}\,,\\
\Im (\hat{\Psi}(z)) &= 2 \ta_0 \sin(\delta) \abs{z}^{\gamma-1}\,,\\
\abs{\hat{\Psi}(z)} &= \abs{z}^{\gamma-1} \sqrt{(1 + 2 \ta_0 \cos(\delta))^2+4 \ta_0^2 \sin(\delta)^2}\;.
\end{aligned}
\end{equation}
Representing $\hat{\Psi}$ in polar coordinates,
\begin{equation*}
\hat{\Psi}(z) = \abs{\hat{\Psi}(z)} \ef{\i \theta(z)}\,,
\end{equation*}
we get
\begin{equation}
\label{theta}
\begin{aligned}
     \sqrt{\hat{\Psi}(z)} &= \sqrt{\abs{\hat{\Psi}(z)}} \ef{\i \theta(z)/2}
\qquad\quad \text{ with }\\
      \theta (z) &= \arctan (\Im (\hat{\Psi}(z))/ \Re (\hat{\Psi}(z))))  \qquad \theta\in (-\pi,\pi) \;.
\end{aligned}
\end{equation}
Note, that $\sqrt{\hat{\Psi}(z)}$ is the complex root with non-negative real part, which meets the general
assumption of the paper.
Moreover, we have
\begin{equation}
\label{szz}
\alpha_{sz}(z) = \Re (\alpha_{sz}^*(z))
          = \eta \Re (\Psi(z)) + \xi \Im (\Psi (z))\;.
\end{equation}

First, we prove that $\Re (\Psi(z)) \geq 0$: We use the elementary inequality
\begin{equation*}
\cos (\theta(z)) \leq \cos^2 (\theta(z)/2)\,,
\end{equation*}
and $\cos (\theta(z)/2)\geq 0$ which imply that
\begin{equation}
\label{eq:real}
\begin{aligned}
    \Re \left(\sqrt{\hat{\Psi}(z)}\right)
      &= \sqrt{\abs{\hat{\Psi}(z)}} \cos (\theta(z)/2) \\
      &=  \sqrt{\abs{\hat{\Psi}(z)}\, \cos^2 (\theta(z)/2) }\\
      &\geq  \sqrt{\Re (\hat{\Psi}(z))}\;.
\end{aligned}
\end{equation}

\begin{itemize}
\item Now, let $\gamma\in (0,1)$.
Since $\Re (\hat{\Psi}(z)) \geq 1$ for $\gamma\in (0,1)$, it follows that for all $z \in \C_0$
\begin{equation*}
    \eta \Re (\Psi(z))
           = \eta \Re \left(\sqrt{\hat{\Psi}(z)}\right) - \eta
          \geq  \eta \sqrt{\Re (\hat{\Psi}(z))} - \eta
\geq 0\,.
\end{equation*}
Thus $\eta \Re (\Psi(z)) \geq 0$.

Now we show that $z \to \xi \Im (\Psi (z))$ is uniformly bounded from below by $0$ in $\C_0$. Thus according to (\ref{szz})
$\alpha_{sz}$ is uniformly bounded from below, and thus from Corollary \ref{co:lion}, it follows that $t \to K(\vx,t)$ is causal.

Using the definition of $\theta$, (\ref{theta}), and the facts that
$\delta \in [0,(1-\gamma)\pi/2]$ for $\phi \in [0,\pi/2]$ and $\delta \in [(\gamma-1)\pi/2,0)$ for $\phi \in (\pi/2,\pi]$
it follows from the monotonicity of $\tan$ on $(-\pi,\pi)$ that
\begin{equation*}
\theta(z) = \arctan\left(\frac{2 \ta_0 \sin(\delta) \abs{z}^{\gamma-1}}{1 + 2 \ta_0 \cos(\delta) \abs{z}^{\gamma-1}}\right) \in
\left\{ \begin{array}{rcl}
~[0,\delta] & \text{ for all } & \phi \in [0,\pi/2]\,,\\
~[\delta,0) & \text{ for all } & \phi \in (\pi/2,\pi]
\end{array}
\right.
\end{equation*}
Now, noting that $\sgn (\xi) = \sgn \left(\sin (\theta(z)/2)\right)$ it follows that
\begin{equation}
\label{eq:imaginary}
\xi \Im (\Psi(z)) = \xi \sqrt{\abs{\hat{\Psi}(z)}} \sin (\theta(z)/2) \geq 0\,.
\end{equation}
Thus the assertion follows from Corollary \ref{co:lion}.

\item Assume that $\gamma  > 1$. Let $z=\xi+\i\eta$ with $\eta = 0$. 
Since the square root in~\req{delta} is such that $\Re(\Psi(z))>0$ for $z\in\C_0$, property~(\ref{propsqrt}) 
in the Appendix implies $\xi \Im (\Psi(z))>0$ for $z=\xi$ and hence $\alpha_{sz}(z=\xi) \geq 0$. 
Thus (\ref{estimate:rapidly}) holds and we can apply Corollary \ref{co:lion}.
\begin{itemize}
\item Let $\gamma\in (1,3)\cup (5,7)\cup\cdots$, which implies that $\cos\left(\gamma\,\pi/2\right)<0$, and consequently, $\ta_0 < 0$.
For sufficiently large $n$ the elements of the sequence $\set{z_n:=\i n}$ satisfy
\begin{equation*}
\alpha_{sz}(z_n) = n\,\Re\left(\sqrt{1 - 2\abs{\ta_0}\, n^{\gamma-1}}-1\right) \leq -n\,,
\end{equation*}
which shows that (\ref{estimate:polynomial_negativ}) holds with $\mu=1$, $C_1=1/2$ and $C_2=0$, and hence the assertion follows
from Corollary \ref{co:lion}.

\item Let $\gamma\in (3,5)\cup (7,9)\cup\cdots$, which implies that $\cos\left(\gamma\,\pi/2\right)>0$, and consequently, $\ta_0>0$.
Now, let $z_n:=n \ef{\i \phi}$ with
$$
    \phi:=  \frac{\pi}{\gamma-1} + \frac{\pi}{2}\,.
$$
Since $\gamma>3$, we have $\phi\in (\pi/2,\pi)$, and therefor $\Re(z_n)=n\,\cos(\phi)<0$ and $\Im(z_n)=n\,\sin(\phi)>0$.
Moreover,
$$
       \hat \Psi(z_n) = 1-2\,\abs{\ta_0}\,n^{\gamma-1}\,
$$
and thus for sufficiently large $n$ we have $\Re(\Psi(z_n))=-1$ and \\$\Im\left(\sqrt{\hat\Psi(z_n)}\right)\geq C_1\,n^{(\gamma-1)/2}$ for
some constant $C_1>0$. Hence it follows that
\begin{equation*}
\begin{aligned}
\alpha_{sz}(z_n)
         &=  -\Im(z_n) + \Re(z_n)\, \Im(\sqrt{\hat\Psi(z_n)})\\
         &= - n\abs{\sin(\phi)} - n\abs{\cos(\phi)}\, C_1\,n^{(\gamma-1)/2}
\end{aligned}
\end{equation*}
and therefore (\ref{estimate:polynomial_negativ}) holds.
Thus from Corollary \ref{co:lion} the assertion follows.
\end{itemize}
\end{itemize}
\end{proof}

\subsection*{Thermo-Viscous Attenuation Law}

\begin{theorem}
Let $c_0,\tau_0>0$ and let $\alpha_{tv}^*$ as defined in~(\ref{alpha*th}).
Then the kernel function $K$ violates causality.
\end{theorem}

\begin{proof}
Since the function $z\in \C_0 \to 1-\i \tau_0\,z$ does not vanish, the function
\begin{equation*}
z \in \C_0 \to \alpha_{tv}^*(z)
     =  \frac{-\i z}{c_0\,\sqrt{1-\i \tau_0 z}} +  \frac{\i z}{c_0}
\end{equation*}
is the holomorphic extension of $\omega \in \R \to \alpha_{tv}^*(\omega)$.
That (\ref{estimate:rapidly}) holds follows from the identity~(\ref{alpha*threal}).
For the sequence $\set{z_n }_{n\in\N} := \set{\i n}_{n\in\N}$ we get for $n$ sufficiently large
$$
   \alpha_{tv}^*(z_n) = \frac{n}{c_0}\,\left[ \frac{1}{\sqrt{1+\tau_0\,n}} -1 \right] \leq - \frac{1}{2\,c_0} n
$$
Thus (\ref{estimate:polynomial_negativ}) holds.

The second part of Corollary \ref{co:lion} implies that $K$ is not causal.
\end{proof}

\subsection*{Model of Nachman, Smith and Waag}

\begin{theorem}
Let $\alpha_{nsw}^*$ as in~(\ref{alpha*Nachman+}). If
\begin{equation}\label{asstau}
  \tilde \tau_m <  \tau_m  \qquad \mbox{ for all} \qquad m\in\{1,\,\ldots\,,N\}\,,
\end{equation}
then the kernel function $K$ is causal.
\end{theorem}

\begin{proof}
Since for all $z \in \C_0$ and all $1 \leq m \leq N$, $1-\i \tau_m\,z$ does not vanish,
\begin{equation*}
\begin{aligned}
z \in \C_0 \to  \alpha_{nsw}^*(z)
     =  \frac{-\i z}{c_0}\,\left[
             \frac{c_0}{\tilde c_0}\, \sqrt{ \frac{1}{N}\, \sum_{m=1}^N
          \frac{1-\i\,\tilde\tau_m\,z}{1-\i \tau_m\,z} }
            -1 \right]
\end{aligned}
\end{equation*}
is the holomorphic extension of $\omega \in \R \to \alpha_{nsw}^*(\omega)$.

We use a similar notation as in Proposition \ref{theo:szabo}.
\begin{equation*}
z = \xi + \i \eta = \abs{z} \ef{\i \phi} \in \C_0\,,
\end{equation*}
with some $\phi \in [0,\pi]$.
\begin{equation*}
\alpha_{nsw}^*(z) = \frac{\abs{z}}{c_0} \ef{\i (\phi-\pi/2)} \,(\Psi(z)-1)\,,
\end{equation*}
where
\begin{equation*}
\Psi(z) = \sqrt{\sum_{m=1}^N \hat{\Psi}_m(z)}\, \quad\mbox{ with }\quad
  \hat{\Psi}_m(z) :=  \frac{1}{N}\, \frac{c_0^2}{\tilde c_0^2}\,
          \frac{1-\i\,\tilde\tau_m\,z}{1-\i \tau_m\,z}  \,,
\end{equation*}

In the following we show that for all $z\in \C_0$
\begin{equation}\label{z2A1z1A2}
c_0 \alpha_{nsw}(z) = \eta (\Re (\Psi)(z)-1) + \xi \Im (\Psi)(z) >0\,,
\end{equation}
which means that (\ref{estimate:polynomial}) holds. Then, according to Corollary \ref{co:lion} the
function $t \to K(\vx,t)$ is causal.

As in the proof of Proposition \ref{theo:szabo} we prove $\eta (\Re (\Psi)(z)-1)>0$ and $\xi \Im (\Psi)(z) >0$.
\begin{itemize}
\item Taking into account~(\ref{proptildec0tau0}) we define
\begin{equation}
\label{s}
   s := \frac{1}{N}\, \frac{c_0^2}{\tilde c_0^2}
             = \left(\sum_{m=1}^N \frac{\tilde\tau_m}{\tau_m}\right)^{-1}\,.
\end{equation}
Using this notation, we get
\begin{equation*}
\begin{aligned}
  \hat{\Psi}_m(z) &=
          s\,\frac{(1+\tilde\tau_m\,\tau_m\,\abs{z}^2) +\i\,(\tau_m\,\bar{z} -\tilde\tau_m\,z)}
           {|1-\i \,\tau_m\,z|^2}  \\
       &=   s\,\frac{\tilde\tau_m}{\tau_m}\,
      \frac{(\tau_m/\tilde \tau_m + \tau_m^2\,\abs{z}^2 + (\tau_m^2/\tilde \tau_m+\tau_m)\,\eta)
                   +\i\,(\tau_m^2/\tilde \tau_m -\tau_m)\,\xi}
           {1+\tau_m^2\,\abs{z}^2 + 2\,\tau_m\,\eta} \,,
\end{aligned}
\end{equation*}
Because $\tau_m/\tilde \tau_m>1$, by assumption (\ref{asstau}), it follows
that for all $z \in \C_0$
$$ \Re(\hat{\Psi}_m(z)) > s\,\frac{\tilde\tau_m}{\tau_m}\;.$$
Consequently, by using the definition of $s$, (\ref{s}), it follows that
$$
    \Re\left(\sum_{m=1}^N \hat{\Psi}_m(z)\right) > s\, \sum_{m=1}^N\frac{\tilde\tau_m}{\tau_m} = 1.
$$
Now, using (\ref{eq:real}) it follows 
\begin{equation*}
    \Re(\Psi(z))
      =    \Re\left(\sqrt{\sum_{m=1}^N \Psi_m(z)}\right)
     \geq  \sqrt{\Re\left(\sum_{m=1}^N \Psi_m(z)\right)} > 1
\end{equation*}
and consequently $\eta (\Re (\Psi(z))-1)>0$.

\item We have
\begin{equation*}
\begin{aligned}
  \Im(\hat{\Psi}_m(z)) =  s\,\frac{\tilde\tau_m}{\tau_m}\,
      \frac{(\tau_m^2/\tilde \tau_m -\tau_m)\,\xi}
           {1+\tau_m^2\,\abs{z}^2 + 2\,\tau_m\,\eta}
\end{aligned}
\end{equation*}
together with the assumption (\ref{asstau}), which state that $\tau_m/\tilde \tau_m>1$,
it follows that $\sgn(\Im(\hat{\Psi}_m(z)))=\sgn(\xi)$. 
According to our assumption, we take that complex root, such that
the real part of the argument is non-negative which together with property~(\ref{propsqrt})
in the Appendix implies
$$
    \sgn\left( \Im\left(\sqrt{\sum_{m=1}^N \hat{\Psi}_m(z)} \right) \right)
     = \sgn\left( \Im\left( \sum_{m=1}^N \hat{\Psi}_m(z)  \right) \right) \,.
$$ 
Therefore,
\begin{equation*}
\begin{aligned}
\sgn (\Im(\Psi(z))) &= \sgn\left(\Im\left(\sqrt{\sum_{m=1}^N \hat{\Psi}_m(z)}\}\right)\right)\\
               &= \sgn\left(\Im\left(\sum_{m=1}^N \hat{\Psi}_m(z)\right)\right)
               = \sgn (\xi)\;.
\end{aligned}
\end{equation*}
This shows the assertion.
\end{itemize}
\end{proof}

\subsection*{Our Model}

\begin{theorem} \label{thcaus01}
For $\alpha_0,\,\tau_0>0$ and $\gamma\in (1,2]$ let $\alpha_{ksb}^*$ be defined as in~(\ref{eq:powlaw2}).
Then $K$, as defined in (\ref{31}), is causal.
\end{theorem}
\begin{proof}
The function $\hat{z}\to 1+(-\i \tau_0 \hat{z})$ does not vanish in $\C_0$. Thus
the holomorphic extension of $\omega \to \alpha_{ksb}^*(\omega)$ is given by
\begin{equation*}
       \hat{z} \in \C_0 \to \alpha_{ksb}^*(\hat{z})
          = \frac{\alpha_0(-\i \hat{z})}{c_0\sqrt{1+(-\i \tau_0 \hat{z})^{\gamma-1}}}\;.
\end{equation*}
In the following let $\hat{z} \in \C_0$. For proving (\ref{estimate:polynomial}) we make a variable transformation
\begin{equation*}
  \alpha_{ksb}^*(\hat{z})
    =  \frac{-\i \tau_0 \hat{z}}{\sqrt{1+(-\i \tau_0 \hat{z})^{\gamma-1}}}
    =  \frac{\alpha_0}{\tau_0 c_0} \frac{-\i z}{\sqrt{1+(-\i z^{\gamma-1})}}\,,
\end{equation*}
and define
$$
       \Psi(z) = \frac{1}{\sqrt{\hat{\Psi}(z)}} \quad\text{ and }\quad
       \hat{\Psi}(z)= 1+(-\i \,z)^{\gamma-1}\;.
$$
Then, with this notation, in order to prove causality of $K$, it suffices to prove that for all $z \in \C_0$
\begin{equation}
\label{positive}
\frac{\tau_0 c_0}{\alpha_0}\alpha_{ksb}(\hat{z}) = \eta \Re (\Psi(z)) + \xi \Im (\Psi(z)) \geq 0\;.
\end{equation}
As in the proof of Proposition \ref{theo:szabo} we show that both terms
$\eta \Re (\Psi(z))$ and $\xi \Im (\Psi(z))$ are non-negative, and then from Corollary \ref{co:lion} the
assertion follows.

In order to prove (\ref{positive}) we note that the function $\hat{\Psi}$ here is the same as in (\ref{eq:delta}) in the proof
of Proposition \ref{theo:szabo} when $\ta_0$ is set to $1/2$. Thus we can already rely on the series of manipulations for $\hat{\Psi}$
developed in the proof of Proposition \ref{theo:szabo}.
\begin{itemize}

\item Since $\eta \geq 0$ it suffices to show that $\Re (\Psi(z)) \geq 0$. We note that for a complex number $a+\i b$
\begin{equation*}
        \Re \left(\frac{1}{a+\i b}\right)
           = \Re \left(\frac{a- \i b}{a^2 + b^2}\right)
           = \frac{1}{a^2 + b^2} \Re(a+\i b)\;.
\end{equation*}
Taking into account the definition of $\Psi$ it therefore suffices to show that $\Re \left(\sqrt{\hat{\Psi}(z)}\right) \geq 0$ in
$\C_0$. Since $\Re \left(\hat{\Psi}(z)\right) \geq 0$ in $\C_0$ for $\gamma\in (1,2]$, it follows that
$\Re \left(\sqrt{\hat{\Psi}(z)}\right) \geq 0$ in $\C_0$.

\item Now, using that
\begin{equation*}
\Im \left(\frac{1}{a+\i b}\right) = \Im \left(\frac{a- \i b}{a^2 + b^2}\right) = -\frac{1}{a^2 + b^2} \Im(a+\i b)\,,
\end{equation*}
it suffices to show that $-\xi \Im (\sqrt{\hat{\Psi}(z)}) \geq 0$ for proving that $\xi \Im (\Psi(z)) \geq 0$.
The proof is along the lines as the analogous part in Proposition \ref{theo:szabo} by taking into account that
here $\gamma \in (1,2)$ (in Proposition \ref{theo:szabo} $\gamma \in (0,1)$). In this case we have now that
sign of $\delta$ is exactly opposite as in the proof of Proposition \ref{theo:szabo}, which in turn gives that
$\Im (\Psi(z))$ has the opposite sign as well, and consequently $-\xi \Im (\sqrt{\hat{\Psi}(z)}) \geq 0$.
Thus the assertion follows from Corollary \ref{co:lion}.
\end{itemize}
\end{proof}

In experiments it has been discovered that several biological tissues satisfy a frequency
power law (\ref{eq:powlaw1b}) with exponent $\gamma \in (1,2)$ (cf.~\cite{Web00,BurRoiBauPal010}).
However, as it has been shown in Theorem \ref{coro:powlaw1}, such models are not causal.
Our proposed model approximates the frequency power law for small frequencies, which is actually
the range where it has been experimentally validated. So, our proposed model, is valid in the
actual range of experimentally measured data and extrapolates the measured data in a causal way.
Figure \ref{fig:comp} shows a comparison of $\alpha_{pl}$ and $\alpha_{ksb}$ in an experimental
frequency range.

\begin{figure}[htb]
\begin{center}
\includegraphics[width=0.45\textwidth]{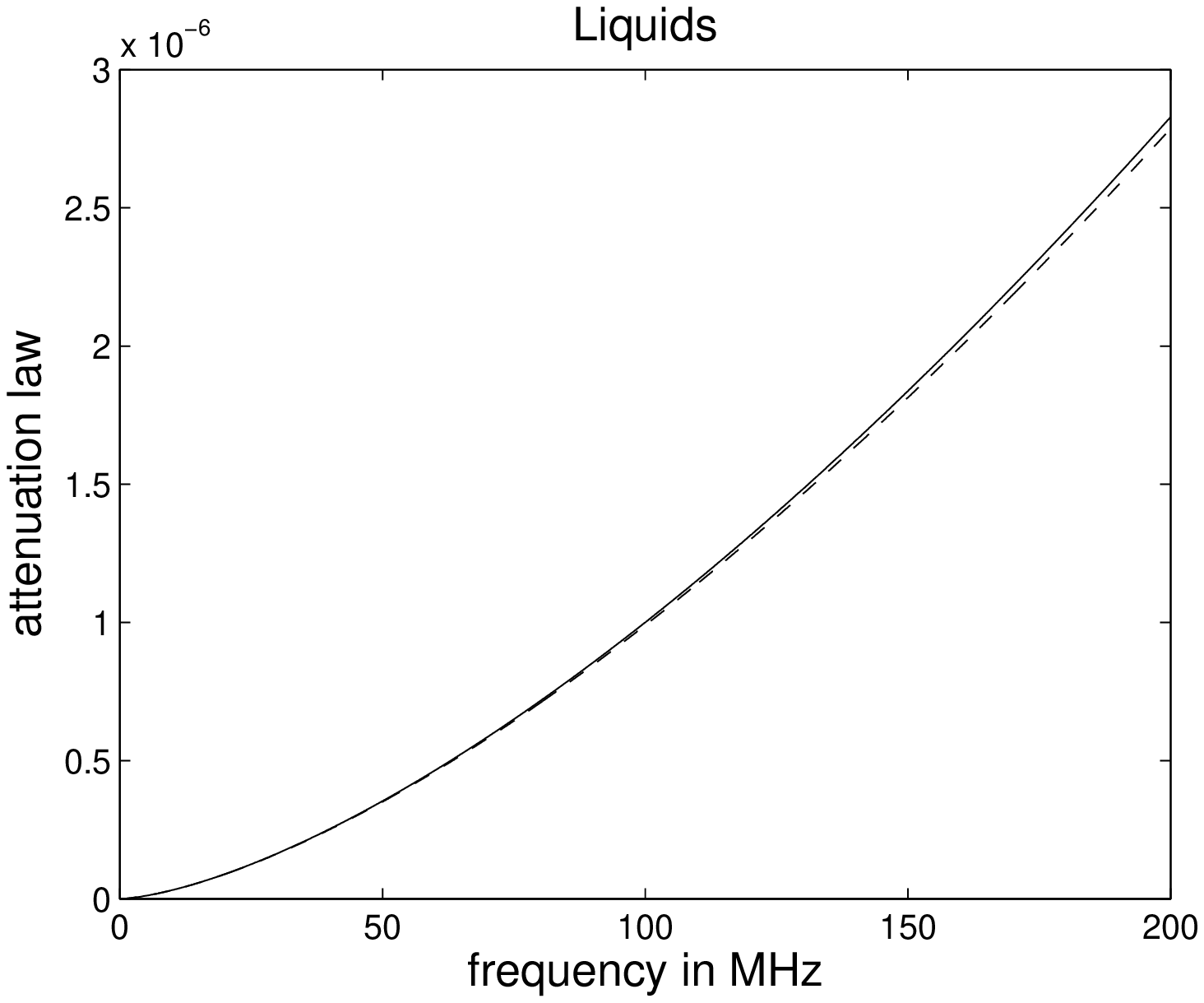}
\includegraphics[width=0.45\textwidth]{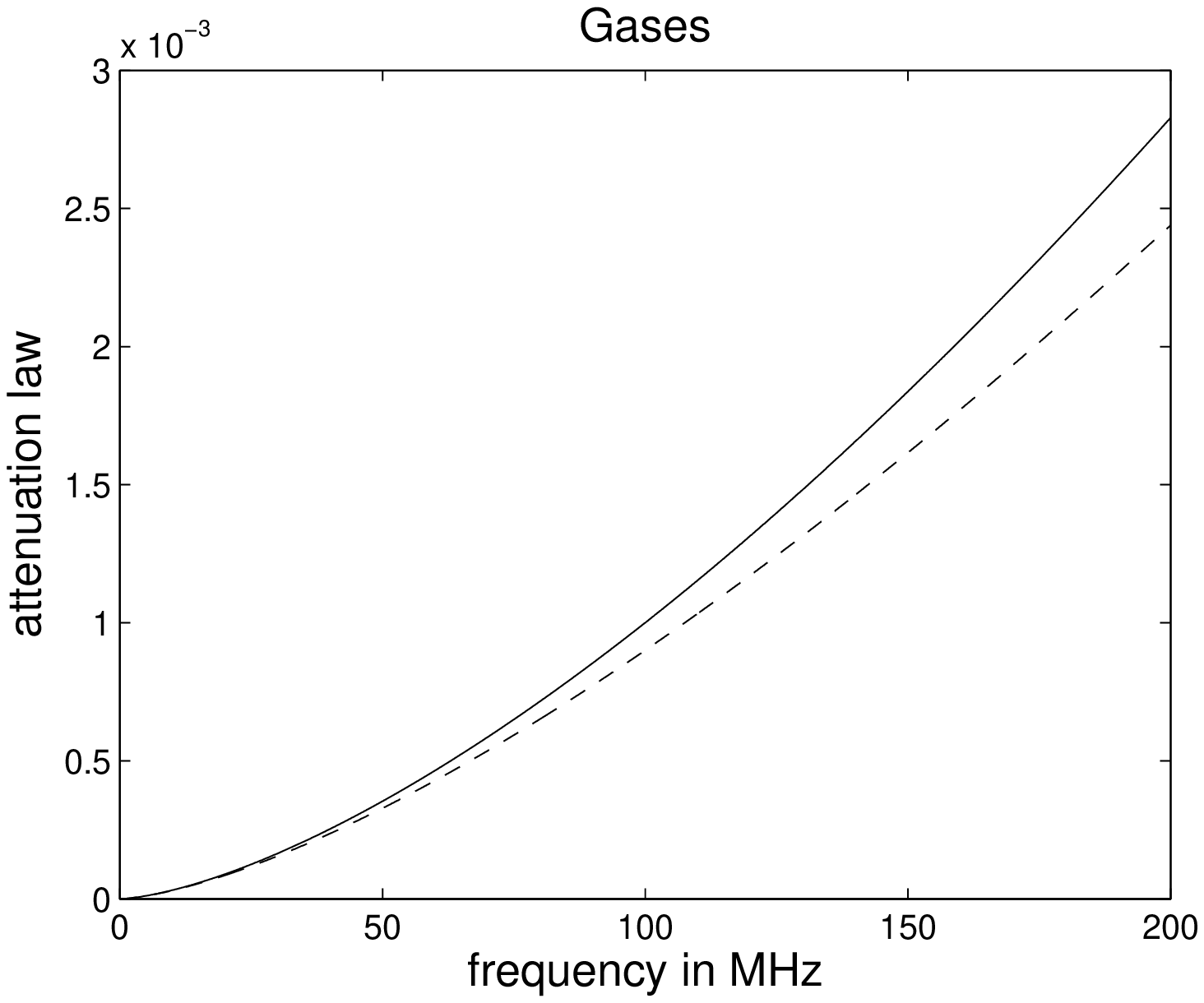}
\end{center}
\caption{For $\gamma=1.5$: Comparison $\omega \to \alpha_{ksb}(\omega)$ (as defined in (\ref{eq:powlaw2})) where
$\alpha_0 := 2c_0\tau_0/\abs{\cos (\frac{\pi}{2} \gamma)}$ (dashed line) and the power law
$\alpha_{pl}(\omega) = \abs{\tau_0\,\omega}^\gamma$  (as defined in (\ref{eq:powlaw1b})).
For liquids: $\tau_0=10^{-6}\,MHz$ (left picture) and for gases: $\tau_0=10^{-4}\,MHz$ (right picture)
(cf.~\cite{KinFreCopSan00}). Experiments for determining the power law coefficient are performed in the range
$0-60\,MHz$ (cf. e.g.~\cite{Sza95}), which is the basis for the range of the represented data.
}\label{fig:comp}
\end{figure}

\subsection*{Model of Greenleaf and Patch}

\begin{proposition}\label{theo:gp}
For $\alpha_0>0$ let $\alpha_{gp}^*$ be defined as in~(\ref{alpha*Patch}) with the specified values
$\gamma\in\set{1,2}$.
Then $K$, as defined in (\ref{31}), is not causal.
\end{proposition}

\begin{proof}
For the two specified models we have $\alpha_{gp1}(\omega)= a_0 \abs{\omega} >0$ and
$\alpha_{gp2}(\omega) = a_0 \omega^2$. The respective holomorphic extensions are given by (cf. Proof of
Theorem~\ref{thcaus02} and Theorem~\ref{coro:powlaw1})
\begin{equation*}
\begin{aligned}
z \in \C_0 \mapsto \hat{\alpha}_{gp1}^*(z) := a_0 z
        + \i \frac{2 a_0}{\pi} z \log \left(\frac{z}{\omega_0}\right) \qquad (\omega_0\neq0, \mbox{fixed})
\end{aligned}
\end{equation*}
and
\begin{equation*}
\begin{aligned}
 z \in \C_0 \mapsto \hat{\alpha}_{gp2}^*(z) := a_0\,(-\i\,z)^2  \,.
\end{aligned}
\end{equation*}
The assertion for $\gamma=1$ follows as in the proof of Theorem~\ref{thcaus02} and the second  assertion
follows from Theorem~\ref{coro:powlaw1} for $\gamma=2$.
\end{proof}

\subsection*{Model of Chen and Holm}

\begin{theorem}\label{th:chenholm}
Let $0<\alpha_1<1$, $\gamma\in (0,2)$ and $\green$ as in (\ref{GreenChenHolm}).
Then there does not exist a constant $c>0$ such that
\begin{equation}\label{CausRef}
  \supp ( \green\left(\cdot,t \right) ) \subseteq B_{c\,t}(\vO)
              \qquad \mbox{ for } \qquad t>0\,,
\end{equation}
i.e. for each $c>0$ the function $t\mapsto \green(\vx,t+\abs{\vx}/c)$ is not causal.
\end{theorem}

\begin{proof}
Let $t>0$ be fixed. Assume that $\vx\mapsto \green(\vx,t)$ has support in $ B_{c\,t}(\vO)$ for some $c>0$.
Then according to the Paley-Wiener-Schwartz Theorem (Cf.~\cite{GasWit99,Hoe03}) the map
$\k\mapsto \ifourier{\green}(\k,t)$ is infinitely differentiable. We show that this is not possible.
According to~(\ref{GreenChenHolm}) and (\ref{AkBk}), we have
\begin{equation*}
\begin{aligned}
  \F_{3D}^{-1}\left\{ \green\right\} (\k,t)
    = \frac{H(t)\,c_0^2}{(2\,\pi)^{3/2}}\,\ef{A(\k)\,t}\, \frac{\sin(B(\k)\,t)}{B(\k)}
\end{aligned}
\end{equation*}
with
\begin{equation*}
    A(\k) := -\alpha_1\,c_0\,\abs{\k}^\gamma\,,
\qquad
    B(\k) := c_0\,\sqrt{\abs{\k}^2-\alpha_1^2\,\abs{\k}^{2\,\gamma}}\,.
\end{equation*}
Since $\gamma \in (0,2)$, the function $\k \mapsto \ef{A(\k)\,t}$
is not infinitely often differentiable at $\k=\vO$ and since the holomorphic function $\frac{\sin(B(\k)\,t)}{B(\k)}$
does not vanish at $\k=\vO$, it follows that $\k\mapsto \F_{3D}^{-1}\{\green\}(\k,t)$ is not infinitely often
differentiable at $\k=\vO$. Consequently, $\vx\mapsto \green(\vx,t)$ cannot have compact support, which
concludes the proof.
\end{proof}

\section{Integro-Differential Equations Describing Attenuation}
\label{sec:integro}

In the following we derive the integro-differential equations for the attenuated pressure $\patt$ for various attenuation laws.
Thereby, we first derive equations which the according attenuated Green functions $\green$ (cf. (\ref{eq:FGG0})) are satisfying, and
then, by convolution, we derive the equations for $\patt$. The integro-differential equations are general in the sense, that they
apply to arbitrary source terms $f$, and in particular to the source term $f$ (\ref{eq:source}) of the forward problem of
photoacoustic imaging with attenuated waves.

For this purpose, we rewrite $\nabla^2 \fourier{\green}$ by using its definition~(\ref{eq:FGG0}), i.e. $\green = K *_t \green_0 $,
and the product differentiation rule, which gives
\begin{equation}\label{eq:laplG}
\begin{aligned}
&\frac{1}{\sqrt{2\,\pi}}\,\nabla^2 \fourier{\green}\\
= & \nabla^2 \fourier{K} \cdot \fourier{\green_0} +
2 \nabla \fourier{K} \cdot \nabla \fourier{\green_0}
+ \fourier{K} \cdot \nabla^2 \fourier{\green_0}.
\end{aligned}
\end{equation}
To evaluate the expression on the right hand side, we calculate $\nabla \fourier{K}$ and $\nabla^2 \fourier{K}$.
From~(\ref{eq:kernel}), it follows that
\begin{equation}\label{eq:defg2b}
\begin{aligned}
  \nabla \fourier{K} = -\beta^*{}' \cdot \fourier{K} \cdot \sgn \,,
\end{aligned}
\end{equation}
where $\beta^*{}'$ denotes the derivative of $\beta^*(r,\omega)$ (cf. (\ref{eq:kernel})) with respect to $r$.
This together with the formula (\ref{eq:der_sgn}) in the Appendix implies that
\begin{equation}\label{eq:defg2c}
\begin{aligned}
& \nabla^2 \fourier{K}\\
         =&- \nabla \cdot \left( \beta^*{}' \cdot \fourier{K} \cdot \sgn \right)\\
         =& - (\nabla \cdot \sgn) \cdot \beta^*{}' \cdot \fourier{K}
            - (\sgn \cdot \nabla \beta^*{}') \cdot \fourier{K}
            - (\sgn \cdot \nabla \fourier{K}) \cdot \beta^*{}' \\
         =& \left[ -\frac{2}{\abs{\vx}} \cdot\beta^*{}'
                   - \beta^*{}''
                   +  \left(\beta^*{}'\right)^2\right] \cdot \fourier{K}.
\end{aligned}
\end{equation}
Inserting \req{defg2b} and \req{defg2c} into \req{laplG} and using again the identity
$\green = K *_t \green_0$, shows that
\begin{equation}\label{eq:laplG2}
 \begin{aligned}
  & \frac{1}{\sqrt{2\,\pi}}\,\nabla^2 \fourier{\green} \\
  &\qquad = \frac{1}{\sqrt{2\,\pi}}\,\left[ -\frac{2}{\abs{\vx}} \cdot \beta^*{}' - \beta^*{}''
      +  \left(\beta^*{}'\right)^2\right]\cdot\fourier{\green}\\
   &\qquad\quad -2 \beta^*{}' \cdot \fourier{K} \cdot (\sgn \cdot \nabla \fourier{\green_0})
      + \fourier{K} \cdot \nabla^2 \fourier{\green_0}.
 \end{aligned}
\end{equation}
From this identity, together with the two following properties of $\green_0$,
\begin{equation}
\label{eq:g1a}
\nabla\fourier{\green_0}
= \left[\frac{\i \omega}{c_0} -\frac{1}{\abs{\vx}}\right] \cdot \fourier{\green_0} \cdot \sgn\,,
\end{equation}
and
\begin{equation}
\label{eq:g1b}
  \nabla^2 \fourier{\green_0} + \frac{\omega^2}{c_0^2}\fourier{\green_0}
      = -\frac{1}{\sqrt{2\pi}}\delta_\vx\,,
\end{equation}
it follows that
\begin{equation}\label{eq:laplG3}
 \begin{aligned}
   \nabla^2 \fourier{\green}
      =& \left[ -\frac{2}{\abs{\vx}} \cdot \beta^*{}'
                   - \beta^*{}''
                   +  \left(\beta^*{}'\right)^2\right] \cdot \fourier{\green}\\
       & - 2 \left[\frac{\i \omega}{c_0} -\frac{1}{\abs{\vx}}\right] \cdot
                   \beta^*{}' \cdot \fourier{\green}\\
       &-\frac{\omega^2}{c_0^2} \cdot \fourier{\green}\\
       &-\fourier{K} \cdot \delta_\vx\;.
 \end{aligned}
\end{equation}
Inserting the identity $\fourier{K}(\vx,\omega) \cdot \delta_\vx = \fourier{K}({\bf 0},\omega) \cdot \delta_\vx$
in \req{laplG3} gives the \emph{Helmholtz equation}
\begin{equation}\label{eq:helmholtz}
\begin{aligned}
~& \nabla^2 \fourier{\green}
       -\left[\beta^*{}' + \frac{(-\i \, \omega)}{c_0}\right]^2 \cdot\fourier{\green}\\
 =& -\beta^*{}'' \cdot \fourier{\green}
           - \fourier{K}({\bf 0},\cdot)\cdot \delta_\vx\\
 =& -\beta^*{}'' \cdot \fourier{\green}
           -\frac{1}{\sqrt{2\pi}} \ef{-\beta^*({\bf 0},\omega)} \cdot \delta_\vx\,.
\end{aligned}
\end{equation}
To reformulate \req{helmholtz} in space--time coordinates, we introduce two convolution operators:
\index{$K_*$} \index{$D_*$}
\begin{equation}\label{eq:defD*D*'}
   D_*f := K_* *_t f \ltext{and}
   D_*'f :=K_*' *_t f ,
\end{equation}
where the kernels $K_*$ and $K_*'$ are given by
\begin{equation}
\label{eq:defKK'}
   K_*:=K_*(\vx,t) := K_*(\abs{\vx},t) \stext{and}
   K_*(r,t)  := \frac{1}{\sqrt{2\pi}} \ifourier{\beta^*{}'}(r,t)
\end{equation}
and
\begin{equation}
\label{eq:defKK'2}
  K_*':=K_*'(\vx,t) := K_*'(\abs{\vx},t)
\stext{and}
    K_*'(r,t)  = \frac{1}{\sqrt{2\pi}} \ifourier{\beta^*{}''}(r,t).
\end{equation}
Using these operators and applying the inverse Fourier transform to \req{helmholtz} gives
\begin{equation}\label{eq:waveeq+}
  \nabla^2 \green
       -\left[D_*
         + \frac{1}{c_0} \frac{\partial}{\partial t}\right]^2 \green
                  =  -D_*' \green
                     - K({\bf 0},\cdot) \delta_\vx\;.
\end{equation}
In the case that $\beta^*(\abs{\vx},\omega) = \alpha^*(\omega)\abs{\vx}$ is of standard form (\ref{eq:attenuation_coefficient}),
it follows that
\begin{equation}
\label{eq:defKK'+}
\begin{aligned}
K_*(t) = \frac{1}{\sqrt{2\pi}} \ifourier{\alpha^*}(t) \qquad\text{ and } \qquad
K_*'\equiv 0\;.
\end{aligned}
\end{equation}

For a general source term $f$, we denote the attenuated wave by $\patt$. \index{$\patt$} That is
\begin{equation*}
\patt := \patt (\vx,t) = \green *_{\vx,t} f =: \A f\,,
\end{equation*}
where $\A$ is the convolution operator according to the Green function $\green$.
This then shows that $\patt$ satisfies the integro-differential equation
\begin{equation}
\label{eq:waveeq+2}
\boxed{
\nabla^2 \patt -\frac{1}{c_0^2} \frac{\partial^2 \patt }{\partial t^2} =  - \A_s f\;,}
\end{equation}
where $\A_s$ denotes the space--time convolution operator with kernel
\begin{equation}\label{eq:defKs}
\begin{aligned}
   K_s:=K_s(\vx,t) :=  -(\B \green)(\vx,t)
           +(D_*' \green)(\vx,t)
          + K({\bf 0},t) \cdot \delta_\vx(\vx)
\end{aligned}
\end{equation}
and
\begin{equation}\label{eq:B}
    \B:= D_*^2 + \frac{2}{c_0} D_*\frac{\partial}{\partial t}\;.
\end{equation}
Equation~\req{waveeq+2} is called \emph{pressure wave equation with attenuation
coefficient} $\beta^*$. We emphasize that $\beta^*$ determines the operators $D_*$ and $D'_*$
which in turn determine the operator $\A_s$, which in turn determines $\patt$ -
this reveals the dependence of $\patt$ from $\beta_*$.

\begin{remark}\label{rema:alpha*2}
Let $\beta^*(r,\omega)=\alpha^*(\omega)r$ be the standard attenuation model (cf. (\ref{eq:attenuation_law})).
Assuming that the associated kernel $K$ (cf. (\ref{31})) is causal, it follows that
\[
 \abs{\nabla K} = \frac{1}{\sqrt{2\pi}} \abs{\ifourier{\alpha^* \cdot \ef{-\alpha^* \abs{\vx}}}}\;.
\]
Using some sequence $\set{\vx_n}$ satisfying $\vx_n \neq {\bf 0}$ and $\vx_n \to {\bf 0}$ shows that
\[
  \lim_{n \to \infty} \abs{\nabla K}(\vx_n,t) = \frac{1}{\sqrt{2\pi}} \abs{\ifourier{\alpha^*}(t)} \underbrace{=}_{(\ref{eq:defKK'+})} \abs{K_*(t)}\;.
\]
Due to the causality of $K$ the left hand side is zero for $t<0$, and thus $K_*$ is also causal.

Because the convolution of causal distributions is well-defined, the operator $D_*$ is well-defined
on all causal distributions.
Moreover, since $K_*'=0$, it follows that $D_*'\equiv 0$. Using that $K_*$ depends only on $t$ it follows that
\[
       (D_* \green) *_{\vx,t} f =  [K_* *_t \green] *_{\vx,t} f
       =  K_* *_t [\green *_{\vx,t} f] = D_*(\green *_{\vx,t} f).
\]
Convolving each term in \req{waveeq+} with a function $f$, using the previous identity and that $D_*' \equiv 0$, it follows that
\begin{equation} \label{eq:waveeq+3}
  \nabla^2 \patt
       -\left[D_* + \frac{1}{c_0} \frac{\partial}{\partial t}\right]^2  \patt = - f \;
\end{equation}
where
\begin{equation}\label{eq:D*sp}
   D_* \cdot = \frac{1}{\sqrt{2\,\pi}}\,\fourier{\alpha^*}(t) *_t \cdot \,.
\end{equation}
\end{remark}

In the following we derive the common forms of the wave equation models corresponding to the
various attenuation models listed in Section~\ref{sec:attenuation}.

\subsection*{Power Laws}
\begin{itemize}
\item Let $0 < \gamma \not\in \N$ and $0 < \alpha_0$.
We note that the \emph{Riemann-Liouville fractional derivative} with respect to time, denote
by $D_t^\gamma$ (see \cite{KilSriTru06,Pod99}), is defined in the Fourier domain by
\index{Riemann-Liouville fractional derivative, $D_t^\gamma$}
\begin{equation}\label{eq:defDtga}
        \fourier{D_t^{\gamma}f} = (-\i \omega)^\gamma \fourier{f}\,,
\end{equation}
and satisfies
\begin{equation}\label{eq:propDtgamma}
    D_t^{2\gamma}f = D_t^\gamma D_t^{\gamma}f \ltext{and}
    \frac{\partial}{\partial t} D_t^{\gamma}f = D_t^\gamma \frac{\partial}{\partial t} f = D_t^{\gamma+1}f.
\end{equation}
From this together with (\ref{eq:powlaw1}) and (\ref{eq:D*sp}), we infer
\begin{equation}\label{deftildealpha0}
        D_* = \tilde\alpha_0\, D_t^\gamma  \qquad \mbox{ with }\qquad
   \tilde \alpha_0 := \frac{\alpha_0}{\cos(\pi\,\gamma/2)}
\end{equation}
and thus wave equation~(\ref{eq:waveeq+3}) reads as follows
\begin{equation}\label{eq:standwaveeq}
\boxed{
\begin{aligned}
  \nabla^2 \patt
       &-\left[\tilde\alpha_0 \,D_t^\gamma
         + \frac{1}{c_0} \frac{\partial}{\partial t}\right]^2  \patt
      = -f  \;.
\end{aligned}}
\end{equation}

\item Let $\gamma=1$, then for the frequency power law (\ref{alphagamma1}) it follows from the
Fourier transform table I in~\cite{Lig64}
\begin{equation*}
\begin{aligned}
         D_* &= -\frac{4\,\alpha_0}{2\,\pi} \left[\,
                    \frac{H(t)}{t^2} - (\log|\omega_0|)\,\delta'_t
                      \right] *_t \\
             &= -\frac{4\,\alpha_0}{2\,\pi}\,  \frac{H(t)}{t^2} *_t \;\;
                + \frac{4\,\alpha_0}{\sqrt{2\,\pi}}\,(\log|\omega_0|)\,\frac{\partial }{\partial t}\,.
\end{aligned}
\end{equation*}
\end{itemize}

\subsection*{Szabo's Attenuation Law:} Let $0 < \alpha_0$ and $0 < \gamma\not\in\N$.  From (\ref{alpha*szabo}) and
(\ref{eq:defDtga}), we get
$$
      \left[D_* + \frac{1}{c_0} \frac{\partial}{\partial t}\right]^2
            = \frac{1}{c_0^2}\frac{\partial^2 }{\partial t^2}
       + \frac{2\,\tilde\alpha_0}{c_0}\,\frac{\partial }{\partial t}\,D_t^\gamma
$$
and thus wave equation~(\ref{eq:waveeq+3}) reads as follows
\begin{equation}\label{szaboseq}
\begin{aligned}
\boxed{
    \nabla^2 \patt -\frac{1}{c_0^2}\frac{\partial^2 \patt}{\partial t^2}
       - \frac{2\,\tilde\alpha_0}{c_0}\,\frac{\partial }{\partial t}\,D_t^\gamma \patt
            = -f(\vx,t)\,.}
\end{aligned}
\end{equation}

\subsection*{Thermo-Viscous Attenuation Law:} From  (\ref{alpha*th}) we get
$$
   \left(\mbox{Id} + \tau_0\,\frac{\partial}{\partial t}\right)\,
          \left[D_* + \frac{1}{c_0} \frac{\partial}{\partial t}\right]^2
            = \frac{1}{c_0^2} \frac{\partial^2}{\partial t^2}
$$
and thus~(\ref{eq:waveeq+3}) becomes
\begin{equation}\label{thviscwaveeq}
\begin{aligned}
\boxed{
    \left(\mbox{Id} + \tau_0\,\frac{\partial}{\partial t}\right)\,\nabla^2 \patt
            - \frac{1}{c_0^2}\frac{\partial^2 \patt}{\partial t^2}
            = -\left(\mbox{Id} + \tau_0\,\frac{\partial}{\partial t}\right)\,f\,.}
\end{aligned}
\end{equation}
This equation is called the \emph{thermo-viscous wave equation}.

\subsection*{Nachman, Smith and Waag \cite{NacSmiWaa90}:}
We carry out the details only for one relaxation process.
\begin{itemize}
\item $N=1$: From
(\ref{alpha*Nachman+}) we get
$$
   \left(\alpha^*(\omega) + \frac{(-\i\,\omega)}{c_0} \right)^2
          = \frac{(-\i\,\omega)^2}{\tilde c_0^2}\,\frac{1-\i\,\tilde \tau_1\,\omega}{1-\i\,\tau_1\,\omega}
$$
which implies
$$
   \left(\mbox{Id} + \tau_1\,\frac{\partial}{\partial t}\right)\,
          \left[D_* + \frac{1}{c_0} \frac{\partial}{\partial t}\right]^2
            = \left(\mbox{Id} + \tilde\tau_1\,\frac{\partial}{\partial t}\right)\,
                        \frac{1}{\tilde c_0^2} \frac{\partial^2}{\partial t^2}\,.
$$
Thus~(\ref{eq:waveeq+3}) reads as follows
\begin{equation}\label{Nachmanwaveeq}
\boxed{
\begin{aligned}
    \left(\mbox{Id} + \tau_1\,\frac{\partial}{\partial t}\right)\,\nabla^2 \patt&
            - \frac{1}{\tilde c_0^2}\left(\mbox{Id} + \tilde\tau_1\,\frac{\partial}{\partial t}\right)\,
                      \frac{\partial^2 \patt}{\partial t^2}\\
            &= -\left(\mbox{Id} + \tau_1\,\frac{\partial}{\partial t}\right)\,f\,.
\end{aligned}}
\end{equation}
If the term with $\tilde\tau_1=0$ is dropped and $\tilde c_0$ is replaced by $c_0$, then we obtain the thermo-viscous
wave equation~(\ref{thviscwaveeq}).

\item $N>1$: For the general case we refer to equation (26) in~\cite{NacSmiWaa90}.

\end{itemize}

\subsection*{Greenleaf and Patch \cite{PatGre06}:}
\begin{itemize}
\item For $\gamma=2$ the attenuation coefficient equals to
$$
        \alpha^*(\omega) = \alpha_0\,\omega^2 = \tilde\alpha_0\,(-\i \omega)^2\,,
$$
where $\tilde\alpha_0$ is defined as in~(\ref{deftildealpha0})
and thus
$$
   D_*  = -\alpha_0\,\frac{\partial^2 }{\partial t^2} = -\alpha_0\,D_t^ 2\,,
$$
which gives wave equation~(\ref{eq:standwaveeq}) with $\gamma=2$.
\item For $\gamma=1$ we have
$$
       \alpha^*(\omega) = \alpha_0\,|\omega| = \alpha_0\,(-\i \omega)\,\i \sgn(\omega)
$$
and thus
$$
    D_* = \alpha_0\,\mathbf{D}^{-1}
        = -\alpha_0\,\frac{\partial}{\partial t} \H \,,
$$
where $\mathbf{D}^{-1}$ and $\H$ denote the \emph{Riesz fractional differentiation operator}
and the Hilbert transform (cf. Appendix), respectively. Therefore the wave equation reads as
follows
\begin{equation*}
\boxed{
\begin{aligned}
  \nabla^2 \patt
       &-\left[\alpha_0\,\mathbf{D}^{-1}
         + \frac{1}{c_0} \frac{\partial}{\partial t}\right]^2  \patt
      = -f  \;.
\end{aligned}}
\end{equation*}
\end{itemize}

\subsection*{Chen and Holm \cite{CheHolm04}:}
Let $\gamma \in (0,2)$. The Green function defined by~(\ref{GreenChenHolm}) satisfies the Helmholtz equation
\begin{equation}\label{Holmeq02}
\begin{aligned}
~ &\frac{\partial^2 \fourier{\green}}{\partial t^2}(\k,t)
   + 2\,\alpha_1\,c_0\,\abs{\k}^\gamma \,\frac{\partial \fourier{\green}}{\partial t}(\k,t)
   +  c_0^2\,\abs{\k}^2 \,\fourier{\green}(\k,t)\\
= & \frac{c_0^2}{(2\,\pi)^{3/2}} \,\delta(t)
\end{aligned}
\end{equation}
for $t\in\R$ and $\k\in\R^3$. Since the fractional Laplacian for a rotational symmetric function
$f$ and $\gamma\in (0,2)$  is defined by (cf. Definition (2.10.1) in~(\cite{KilSriTru06})
\begin{equation*}
\begin{aligned}
   & \left(-\nabla^2\right)^{\gamma/2} f(\vx)\\
:= &\frac{1}{\sqrt{(2\,\pi)^3}}\,\int_{\R^3} \ef{\vx\cdot\k}\,
                          \left[ \abs{\k}^\gamma\,
     \frac{1}{\sqrt{(2\,\pi)^3}}\,\int_{\R^3} \ef{-\vx\cdot\k}\, f(\vx)\,\d\vx \right]
             \,\d\k \\
= & \F_{3D}\left\{ \abs{\k}^\gamma\,\F_{3D}^{-1}\{ f\}(\k) \right\}(\vx)\,,
\end{aligned}
\end{equation*}
we obtain the following wave equation for $\patt:=\green\ast_{\vx,t}f$
\begin{equation}\label{Holmeq01}
\boxed{
\begin{aligned}
    \nabla^2 \patt -\frac{1}{c_0^2}\frac{\partial^2 \patt}{\partial t^2}
       - \frac{2\,\alpha_1}{c_0}\,\frac{\partial }{\partial t}\,\left(-\nabla^2\right)^{\gamma/2} \patt
            = -f(\vx,t)\;.
\end{aligned}}
\end{equation}
We note that Chen and Holm used instead of $2\,\alpha_1/c_0$ the term $2\,\alpha_1/c_0^{1-\gamma}$.

\subsection*{Our Model~\cite{KowSchBon10}:} From~(\ref{eq:powlaw2}) we get
$$
   \left(\mbox{Id} + \tau_0^{\gamma-1}\,D_t^{\gamma-1}\right)\,
          \left[D_* + \frac{1}{c_0} \frac{\partial}{\partial t}\right]^2
            = \frac{1}{c_0^2} \frac{\partial^2}{\partial t^2} \,
            \left(\alpha_0\,\mbox{Id} + L^{1/2}\right)^2
$$
where the time convolution operator $L^{1/2}$ is the convolution operator with kernel
$$
 l(t):=L^{1/2} (\delta_t)
   = \frac{1}{\sqrt{2\,\pi}}\, \F\left\{ \sqrt{1+(-\i \tau_0\,\omega)^{\gamma-1}} \right\}\;.
$$
Consequently, $L:=(L^{1/2})^2=\mbox{Id} + \tau_0\,D_t^{\gamma-1}$ and (\ref{eq:waveeq+3}) can be rewritten as follows
\begin{equation}\label{causthviscwaveeq}
\boxed{
\begin{aligned}
    \left(\mbox{Id} + \tau_0^{\gamma-1}\,D_t^{\gamma-1}\right)\,\nabla^2 \patt &
            - \frac{1}{c_0^2}\,\left(\alpha_0\,\mbox{Id} + L^{1/2}\right)^2
                       \,\frac{\partial^2 \patt}{\partial t^2}\\
            &= -\left(\mbox{Id} + \tau_0^{\gamma-1}\,D_t^{\gamma-1}\right)\,f\,.
\end{aligned}}
\end{equation}

\section{Pressure Relation}
In this section we derive the relation between $\patt$ and $p_0$ when the source term is of the form (\ref{eq:source}). This chapter is a special instance of
Section \ref{sec:integro}. However, utilizing the special structure of the source term different formulas can be derived.

Attenuation is defined as a multiplicative law (in the frequency domain) relating the
amplitudes of an attenuated and an unattenuated wave initialized by a delta impulse.
Here we are concerned in deriving the convolution relation between the solution $p_0$
of \req{ex:wave3d} (or equivalently of ~\req{ex:ivp3d} and~\req{init_values})
and the attenuated wave function $\patt$, which, according
to (\ref{eq:FGG0}) and (\ref{eq:a_conv}), is given by
\begin{equation}\label{eq:defpatt}
 \patt = \green *_{\vx,t} f  = (K *_t \green_0) *_{\vx,t} f\,,
\end{equation}
with $f$ from (\ref{eq:source}).
Using (\ref{eq:g0}) and the rotational symmetry of $K$, it follows that
\begin{equation*}
\begin{aligned}
  &  \qquad  \qquad \qquad \qquad  \qquad \qquad\patt(\vx,t) \\
       &= \int_\R \int_{\R^3} (K \ast_{t}\green_0)(\vx-\vx',t-t'')
                    \rho(\vx')\,d\vx' \frac{\partial \delta_t}{\partial t}(t'')\,dt''\\
       &= \int_\R \int_{\R^3} \int_\R K(\vx-\vx',t-t'-t'')\green_0(\vx-\vx',t')\,dt'
                    \rho(\vx')\,d\vx' \frac{\partial \delta_t}{\partial t}(t'')\,dt''\\
       &= \int_\R \int_{\R^3} \int_\R \frac{\partial }{\partial t}K(\vx-\vx',t-t'-t'')\green_0(\vx-\vx',t')
                    \rho(\vx') \delta_t (t'')\,dt'\,d\vx' \,dt''\\
    &= \int_{\R^3} \int_\R \frac{\partial }{\partial t} K(\abs{\vx-\vx'},t-t')\,
    \frac{\delta_t(t'-\abs{\vx-\vx'}/c_0)}{4\pi \abs{\vx-\vx'}} \,\rho(\vx')\,\d t'\,\d\vx'\\
    &= \int_{\R^3}\frac{\partial }{\partial t} K(\abs{\vx-\vx'},t-\abs{\vx-\vx'}/c_0)\,\frac{\rho(\vx')}{4\pi \abs{\vx-\vx'}} \,\d\vx'\,.\\
\end{aligned}
\end{equation*}

Using the representation $\vx'-\vx= r' \vec{s}$ with $r' \geq 0$ and $\vec{s} \in S^2$, it follows that
\begin{equation}\label{eq:pattrep01}
\patt(\vx,t) = \frac{1}{4 \pi}\, \int_0^\infty \frac{\partial }{\partial t} K(r',t-r'/c_0) \,r'\,
                                \int_{S^2}  \rho(\vx+r'\vec{s})\,d\vec{s}\,dr'\;.
\end{equation}
Moreover,
\begin{equation*}
\begin{aligned}
p_0(\vx,t)
&= \frac{\partial }{\partial t} \int_{\R^3} \frac{\delta_t(t-\abs{\vx-\vx'}/c_0)}{4\pi \abs{\vx-\vx'}} \rho(\vx')\,\d\vx'\\
&= \frac{\partial }{\partial t} \int_0^\infty r'^2 \int_{S^2}  \frac{\delta_t(t-r'/c_0)}{4\pi r'} \rho(\vx + r' \vec{s})\,d\vec{s}\,dr'\\
&= \frac{\partial }{\partial t} \int_0^\infty r' \frac{\delta_t(t-r'/c_0)}{4\pi} \int_{S^2} \rho(\vx + r' \vec{s})\,d\vec{s}\,dr'\\
&= \frac{\partial }{\partial t} \int_0^\infty \frac{c_0^2 r''}{4\pi} \delta_t(t-r'') \int_{S^2} \rho(\vx + c_0 r'' \vec{s})\,d\vec{s}\,d r''\\
&= \frac{\partial }{\partial t} \left(\frac{c_0^2 t}{4\pi} \int_{S^2} \rho(\vx + (c_0t) \vec{s})\,d\vec{s} \right)\;.
\end{aligned}
\end{equation*}
This gives
\begin{equation}
\label{eq:summ}
\begin{aligned}
  r'\, \int_{S^2}  \rho(\vx+r'\vec{s})\,d \vec{s}
  =  \frac{4\,\pi}{c_0} \int_0^{r'/c_0} p_0(\vx,t')\,\d t'\;.
\end{aligned}
\end{equation}
Now, denoting
\begin{equation*}
F(t,r') := \int_0^{r'}  \frac{\partial}{\partial t} K(r'',t-r''/c_0) \,dr''\,,\qquad
G(r') := \int_0^{r'/c_0} p_0(\vx,r'')\,dr''\,
\end{equation*}
and
\begin{equation*}
F(t,\infty) = \lim_{r' \to \infty} F(t,r')
\end{equation*}
it follows from (\ref{eq:pattrep01}) and (\ref{eq:summ}) and the fact that $p_0(\vx,0) = 0$
that
\begin{equation}
\label{eq:pattrep02a}
\begin{aligned}
\patt(\vx,t) =& \frac{1}{c_0} \int_0^\infty F'(t,r') G(r')\,dr'\\
=& - \frac{1}{c_0} \int_0^\infty F(t,r') \underbrace{G'(r')}_{=p_0(\vx,r'/c_0)/c_0}\,dr' + \frac{1}{c_0} \left. F(t,r') G(r') \right|_{r'=0}^\infty\\
=& \frac{1}{c_0^2} \left( F(t,\infty) \int_0^\infty p_0(\vx,r'/c_0)\,dr'- \int_0^\infty F(t,r') p_0(\vx,r'/c_0)\,dr'\right)\\
=& \frac{1}{c_0} \left( F(t,\infty) \int_0^\infty p_0(\vx,t')\,dt' - \int_0^\infty  F(t,c_0 t') p_0(\vx,t')\,d t'\right)\\
=:& \int_0^\infty \M(t,t') p_0(\vx,t')\,dt'\,,
\end{aligned}
\end{equation}
where
\begin{equation}\label{eq:mos}
\begin{aligned}
  \M(t,t') :&= \frac{1}{c_0} \,(F(t,\infty)-F(t,c_0 t')) \,.
\end{aligned}
\end{equation}
\index{$\M$}
In the following we derive an equivalent representation of $\M$ in terms of the attenuation coefficient,under the assumption
that the attenuation coefficient $\alpha^*$ is such that $K$ is causal.
From \req{kernel}, \req{attenuation_law} and Item~\ref{item:Fdelta} in the Appendix, it follows that
\begin{equation*}
\begin{aligned}
  K(r',t-r'/c_0)
    = \frac{1}{\sqrt{2\,\pi}}\,  \ifourier{ \ef{-\alpha^*(\omega)\,r'+\i\,\frac{\omega}{c_0}\,r'}  }(t) \,.
\end{aligned}
\end{equation*}
which implies
\begin{equation}\label{eq:pattrep02b}
\begin{aligned}
~ & F(t,c_0\,t')\\
    =& \int_0^{c_0\,t'} \frac{\partial }{\partial t} K(r',t-r'/c_0)\,\d r' \\
    =& \frac{1}{\sqrt{2\,\pi}}\,  \ifourier{\frac{-\i\,\omega\,\left[ \ef{-\alpha^*(\omega)\,c_0\,t'+\i\,\omega\,t'} -1 \right]}
                                                 {-\alpha^*(\omega)+\i\,\omega/c_0}\, }(t)  \\
    =& \frac{1}{\sqrt{2\,\pi}}\, \ifourier{\frac{-\i\,\omega\,}{-\alpha^*(\omega)+\i\,\omega/c_0}\, }(t)  \ast_t
          \left[ K(c_0\,t',t-t')  - \delta_t(t)\right]\,.
\end{aligned}
\end{equation}
Since $K$ is causal, it satisfies $K(c_0\,t',t-t')   = 0$ for $t < t'$, and therefore
$$
     F(t,\infty) = \lim_{t'\to \infty} F(t,c_0\,t')
          =  \frac{1}{\sqrt{2\,\pi}}\, \ifourier{\frac{-\i\,\omega\,}{-\alpha^*(\omega)+\i\,\omega/c_0}\, }(t)\,.
$$
Hence \req{pattrep02b} can be written as follows:
\begin{equation*}
\begin{aligned}
  F(t,c_0\,t')
    &= F(t,\infty)  \ast_t
          K(c_0\,t',t-t')  - F(t,\infty)\,.
\end{aligned}
\end{equation*}
and therefore \req{mos} simplifies to
\begin{equation}\label{eq:defM}
\begin{aligned}
   \M(t,t')
      &=   - \frac{1}{c_0}\,F(t,\infty)  \ast_t K(c_0\,t',t-t') \\
      & =\frac{1}{\sqrt{2\,\pi}}\,  \ifourier{\frac{\i\,\omega \,\ef{-\alpha^*(\omega)\,c_0\,t'+\i\,\omega\,t'}}
                                                  {-\alpha^*(\omega)\,c_0+\i\,\omega}\, }(t) \,.
\end{aligned}
\end{equation}
Note that $\M(t,0)= - \frac{1}{c_0}\,F(t,\infty)$.
The following lemma shows that if $K$ is causal
\begin{equation}\label{eq:propM}
         \M(t,t') = 0 \qquad \mbox{ if } \qquad 0< t < t'\,
\end{equation}
and therefore the upper limit of integration in the last term \req{pattrep02a} can be replaced by $t$.
This means that the set of attenuated pressure values
$$\set{\patt(\vx,s) : 0\leq s\leq t}$$ depend only on
the unattenuated pressure values
$$\set{p_0(\vx,s) : 0\leq s\leq t}\;.$$

\begin{lemma}\label{lemm02}
Let $K$ from~\req{kernel} be causal with $\beta^*(r,\omega):=\alpha^*(\omega)\,r$.
Moreover, assume that $\alpha^*(\omega) \neq \i\,\omega/c_0$ for every $\omega\in\R$, and let $\M$ be as defined in~\req{defM}. Then,
\begin{itemize}
\item the function $t \to \M(t,0)$ is causal.
\item For every $t'>0$
\begin{equation}\label{propM1}
    \M(t,t')=0 \text{ for all } t < |t'|\,.
\end{equation}
\end{itemize}
\end{lemma}

\begin{proof}
\begin{itemize}
\item In order to prove causality of $t \to \M(t,0)$ we verify the three assumptions of Theorem~\ref{th:lion} for the tempered distribution
$$
    \fourier{\M}(\omega,0) = \frac{1}{\sqrt{2\,\pi}}\,\frac{\omega }{k(\omega)\,c_0}
 \qquad\mbox{with}\qquad  k(\omega):= \i\,\alpha^*(\omega) + \frac{\omega}{c_0}\,.
$$
Since $K$ is causal, as has been shown in Remark~\ref{rema:alpha*2}, also the function
$$t\mapsto K_*(t)=\frac{1}{\sqrt{2\pi}} \ifourier{\alpha^*}(t)$$ is causal. Now, using Theorem~\ref{th:lion}, it follows that
\begin{enumerate}
\item [1)] $\alpha_*$ is holomorphic in $\mathring{\C}_0$,
\item [2)] $\alpha^*(\xi+\i\,\eta)\to \alpha^*(\xi)$ for $\eta\to 0$ in $\S'$ and
\item [3)] for each $ \epsilon>0$
           there exists a polynomial $P$ such that $\abs{\alpha^*(z)}\leq P(\abs{z})$ for $z\in\C_\epsilon$.
\end{enumerate}
Since $\alpha_*(\omega)\neq \i\,\omega/c_0$ for all $\omega\in \R$ together with
in~\cite[Theorem 2.7]{BelWoh66} it follows that $z \to \alpha_*(z)$ is unique holomorphic extension to $\C_0$ and therefore $z \to \alpha_*(z)$
cannot be identical to $z \to \i\,z/c_0$, the holomorphic extension of $\i\,\omega/c_0$.
$\alpha_*(z)\neq\i\,z/c_0$ for $z\in \C_0$ implies that $k$ has no zeros and hence $z/k(z)$
is holomorphic on $\mathring{\C}_0$. This shows that Item~\ref{it1_Lions} in Theorem~\ref{th:lion} is satisfied for
$t \to \M(t,0)$.

Since $1/(k(\xi+\i\,\eta)\,k(\xi))$ is bounded and $k(\xi+\i\,\eta)\to k(\xi)$ for $\eta\to 0$ in $S'$, it
follows that for all $\psi\in \S$
$$
   \lim_{\eta\to 0} \int_\R \left[ \frac{k(\xi)-k(\xi+\i\,\eta)}{k(\xi+\i\,\eta)\,k(\xi)}\right]\,\psi(\xi)\,\d\xi \to 0
$$
i.e. $1/k(\xi+\i\,\eta)\to 1/k(\xi)$ for $\eta\to 0$ in $S'$. Hence
$\fourier{\M}(\xi+\i\,\eta,0)\to \fourier{\M}(\xi,0)$ for $\eta\to 0$ in $\S'$, which shows Item~\ref{it2_Lions}
in Theorem~\ref{th:lion}.

Since $k(z)$ does not vanish on $\C_0$ and $\abs{k(z)}$ is bounded by a polynomial in $\abs{z}$ for
$z\in\mathring{\C}_0$, it follows that $\abs{z/k(z)}$ is bounded by a polynomial in $\abs{z}$. Hence
Item~\ref{it3_Lions} in Theorem~\ref{th:lion} is satisfied and consequently $t \to \M(t,0)$ is causal.

\item Property~(\ref{propM1}) is satisfied if
\begin{equation}\label{MMM}
\begin{aligned}
    \M(t + \abs{t'},t')=0 \qquad \mbox{ for } \quad t < 0\,.
\end{aligned}
\end{equation}
From~\req{defM} and
$$
K(c_0\,t',t)=\F^{-1}\left\{ \frac{e^{-\alpha^*(\omega)\,c_0\,t'}}{\sqrt{2\,\pi}}\,  \right\}(t)\,,
$$
it follows that
\begin{equation*}
      \M(t + \abs{t'},t')
          = \M(t,0) \ast_t K(c_0\,t',t)\,.
\end{equation*}
Since $t\mapsto\M(t,0)$ and $t\mapsto K(c_0\,t',t)$ are causal, their convolution is also causal (cf. Item~\ref{item:conv}
in the Appendix). This proves property~(\ref{MMM}) and concludes the proof.
\end{itemize}
\end{proof}

\begin{remark}
Assume that the attenuation coefficient is given by
$$
        \alpha^*(\omega) = \i\,\omega/c_0  \qquad \mbox{ for }\qquad \omega\in\R\,  .
$$
Then
$$
   K(r,t) = \frac{1}{\sqrt{2\,\pi}}\,\ifourier{\ef{-\i\,\omega\,r/c_0}}(t)  = \delta(t+r/c_0)
$$
which implies together with~\req{FGG0} and~\req{g0} that
$$
    \green(r,t) = \frac{\delta(t)}{4\,\pi\,r} \,.
$$
But this function does not correspond to the intuition of an attenuated wave, which is manifested by
the convolution equation (\ref{eq:FGG0}), which should give a smooth decay of frequency components over travel
distance. With this Green function $\green$ the input impulse collapses immediately and consequently, in this case,
the assumption $\alpha^*(\omega)\neq\i\,\omega/c_0$ in Theorem \ref{lemm02} reflects physical reality.
\end{remark}

\section{Solution of the Integral Equation}

The inverse problem of photoacoustics with attenuated waves reduces to solving the integral equation (\ref{eq:pattrep02a})
for $p_0$, and to the standard photoacoustical inverse problem, which consists in calculating the initial pressure $\rho$
in the wave equation (\ref{eq:ex:ivp3d}) from measurements of $p_0(\vx,t)$ over time on a manifold surrounding the object
of interest. The standard photoacoustical imaging problem is not
discussed here further, but we focus on the the integral equation (\ref{eq:pattrep02a}).

In the following we investigate the ill--conditionness of the integral equation (\ref{eq:pattrep02a}), where the kernel
$\M$ is given from the attenuation law \req{powlaw2} with $\gamma \in (1,2]$. In this case the
model is causal and the parameter range $\gamma \in (1,2]$ is relevant for biological imaging.

In order to estimate the ill--conditionness of the integral equation (\ref{eq:pattrep02a}) it is rewritten in Fourier domain:
\begin{equation}
\label{eq:Mfrequency}
    \fourier{\patt}(\vx_0,\omega)
         =  \int_\R  \fourier{\M} (\omega,t')\,p_0(\vx_0,t')\,\d t'\;.
\end{equation}
After discretization the ill-conditionness of this equation is reflected by the decay rate of the singular values of
the matrix $\fourier{\M} (\omega,t')$ at certain discrete frequencies and time instances.

We consider simple test examples of attenuation coefficients $\rho$ (as in (\ref{eq:u(x)})), which are characteristic
functions of balls with center at the origin and radiii $R$.
For these examples we investigate the dependence of the ill--conditionedness of (\ref{eq:Mfrequency}) on the radius $R$ and
the location $\vx_0$. For applications in photoacoustic imaging $\vx_0$ would be the location of a detector
outside of the object of interest, to be imaged. Then, by solving the integral equation (\ref{eq:Mfrequency}) $p_0$ can be
calculated, and in turn, the absorption energy $\rho$ can be reconstructed with standard backprojection formulas.
Since
\begin{equation*}
   \supp (p_0(\vx_0,\cdot))
          = \left[(\abs{\vx_0}-R)/c_0,(\abs{\vx_0}+R)/c_0\right]\,,
\end{equation*}
the integral equation~\req{Mfrequency} can be rewritten as
\begin{equation}\label{eq:Mfrequency*}
    \fourier{\patt}(\vx_0,\omega)
         =  \int_{(R_0-R)/c_0}^\infty  \fourier{\M} (\omega,t')\,p_0(\vx_0,t')\,\d t'\,,
\end{equation}
where $R_0 = \abs{\vx}_0$.

In the following we analyze the integral equation (\ref{eq:Mfrequency*}) in terms of the two parameters $R$
and $R_0=\abs{\vx_0}$. This gives a clue on the effect of attenuation in terms of the size
of the object and the distance of the location $\vx_0$ to the simple object. In order to show the effect
of attenuation on the single frequencies, we make a singular value decomposition of the kernel of the
integral equation (\ref{eq:Mfrequency*}).

\begin{example}
\label{exam:distance2}
For small frequencies the attenuation law of castor oil, which behaves very similar to biological soft tissue,
is approximately a power law with exponent $\gamma=1.66$ and $\hat{\alpha}_0 \approx 4\cdot 10^{-2}\,\frac{1}{cm\,(MHz)^\gamma}$, i.e.
$$
   \alpha_{pl} (\omega) \approx  4\cdot 10^{-2}\cdot \omega^{1.66} \cdot cm^{-1} \qquad \mbox{($\omega$ in $MHz$)}.
$$
The sound speed of castor oil is $1490\cdot\frac{m}{s}$ at $25$ degree Celsius. In units of
$cm$ and $MHz$ we have
$$
    c_0\approx 0.15\cdot cm\cdot MHz\,.
$$
Since ~\req{powlaw2} approximates the power law (cf. Figure \ref{fig:comp}) it follows that
$$
  \hat{\alpha}_0\,\abs{\omega}^\gamma \approx
    \frac{\alpha_0\,\sin(\frac{\pi}{2}\,(\gamma-1))}{2\,c_0\,\tau_0}\,\abs{\tau_0\,\omega}^\gamma
$$
and consequently the coefficients of ~\req{powlaw2} satisfy
$$
        \alpha_0 \approx \frac{2\,c_0\,\hat{\alpha}_0}
                        {\tau_0^{(\gamma-1)}\,\sin(\frac{\pi}{2}\,(\gamma-1))}
                 \approx 6\,.
$$
We note that the relaxation time is $\tau\approx 10^{-4}\,\frac{1}{MHz}$ for liquids (cf.~\cite{KinFreCopSan00}).

For the calculation of the singular value decomposition of the discretized kernel of the integral equation
(\ref{eq:Mfrequency*}) we used a frequency range \\
$\omega \in [-80,80] MHz$ and
step size $\Delta \omega=\frac{2\,\pi}{N-1}\, MHz$ with $N=2^9$. The time interval has been set to
$[0,\frac{2\,\pi}{\Delta \omega}] MHz^{-1}$ and a step size $\Delta t=\frac{2\,\pi}{80}\,MHz^{-1}$ was used.

The upper left picture in Fig.~\ref{fig:distance2} visualizes the discretized kernel of
the integral equation (\ref{eq:Mfrequency*}) for $R_0=R$, i.e., when $\vx_0$ is directly on the surface of the
object of interest. The upper right picture shows the singular values of the discretized kernel in a logarithmic scale.
Two properties of the singular values become apparent:
\begin{enumerate}
\item For large indices the decay rate is exponential, which can be seen from the linear decay in the logarithmic scale.
\item Secondly, there is a range of indices, where the singular values do not decay that rapidly. As a consequence, for
      solving the integral equation this means that the Fourier coefficients of $p_0$ according to the first block of singular values can be
      determined in a stable manner.
\end{enumerate}
For increasing distance $L=\abs{\vx_0}-R$ of $\vx_0$ to the object the singular values of the discretization of
the integral equation (\ref{eq:Mfrequency*}) show a drastically more exponentially decay rate for increasing $L$
(see bottom right picture Fig.~\ref{fig:distance2}).
This means that if the object is further away from $\vx_0$ attenuation is more drastically, and solution of the
integral equation is more unstable.
We analyze the dependence of the number of largest singular values from $L$. For this purpose we denote by
$n_{cut}$ the index of the singular value that is about $0.1\%$ of the maximal singular value. For the numerical solution of
(\ref{eq:Mfrequency*}) it means that if we make a truncated singular valued decomposition with only $n_{cut}$ singular values,
the error amplification can be bounded by a factor $1000$. The dependence of $n_{cut}$ on $L$ is shown in the lower left
picture of Fig.~\ref{fig:distance2}. The picture reveals that for increasing distance (from about $2cm$) only about four Fourier
modes of $p_0$ are significant when a maximal error amplification of a factor $1000$ is required. This reveals that in general the
solution of the integral equation (\ref{eq:Mfrequency*}) is significantly ill--posed and worse if the data recording is far away from the
object.
\begin{figure}[htb]
\begin{center}
\includegraphics[angle=90,width=\textwidth]{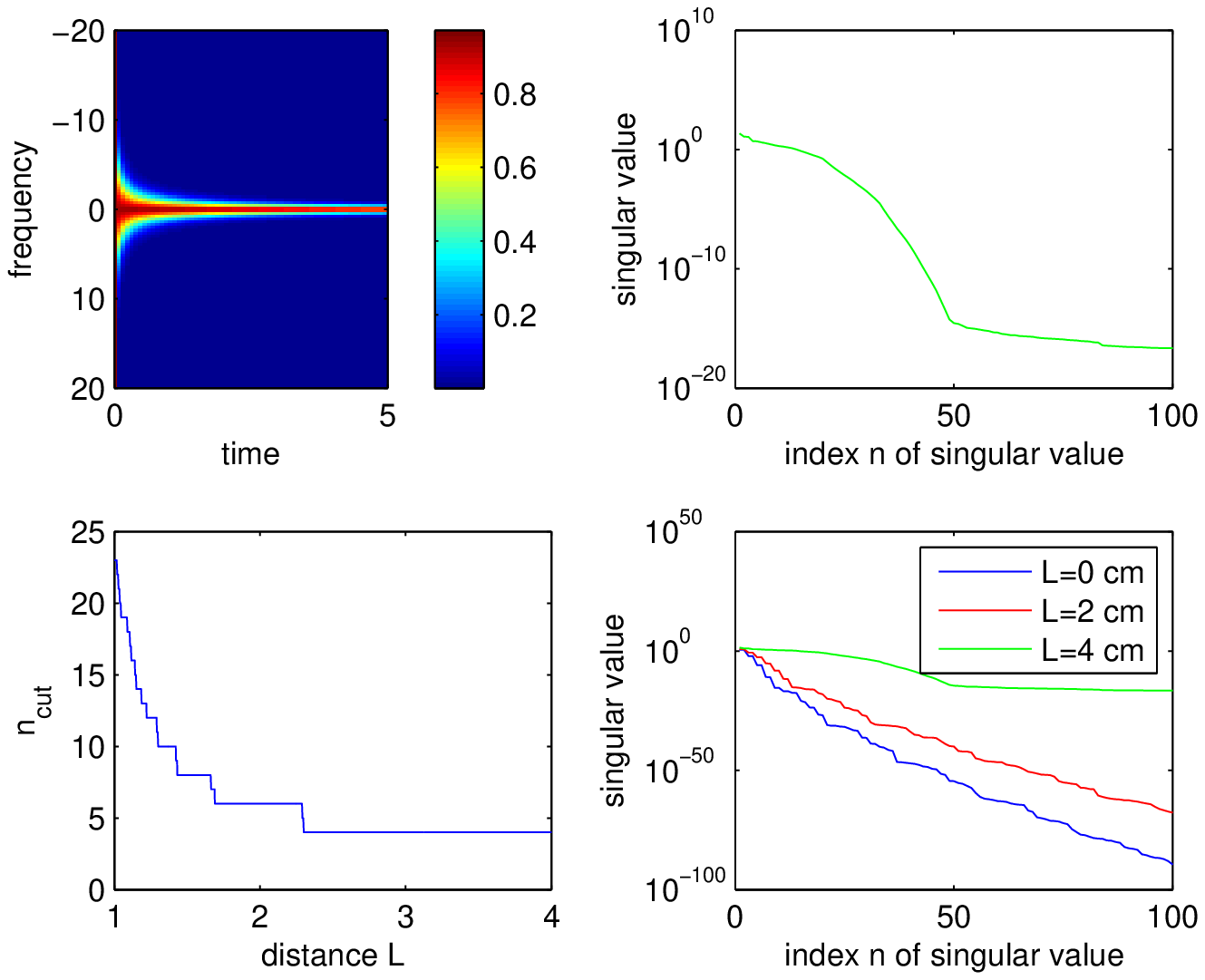}
\end{center}
\caption{Case: $\gamma=1.66$ (castor oil). The upper left and right pictures visualize
the kernel $\fourier{\M}(\omega,t')$ and its singular values for $L:=R_0-R=0cm$.
The lower right and left pictures visualize $\fourier{\M}(\omega,t')$ for the detector
distances $L=0\cdot cm$, $L=2\cdot cm$ and $L=4\cdot cm$ and the respective indices $n_{cut}$ for
which the singular values are about $0.1$ per cent of the maximal singular value.}
\label{fig:distance2}
\end{figure}
\end{example}

\begin{example}\label{exam:distance1}
An analogous numerical example as in Example~\ref{exam:distance2} for the case $\gamma=1.1$
is presented in Fig.~\ref{fig:distance1}. From the lower left picture of Fig.~\ref{fig:distance1},
we see that if the distance is about $2\cdot cm$ from the boundary of the object,
then $17$ singular values are available for the numerical estimation.
\begin{figure}[htb]
\begin{center}
\includegraphics[angle=90,width=\textwidth]{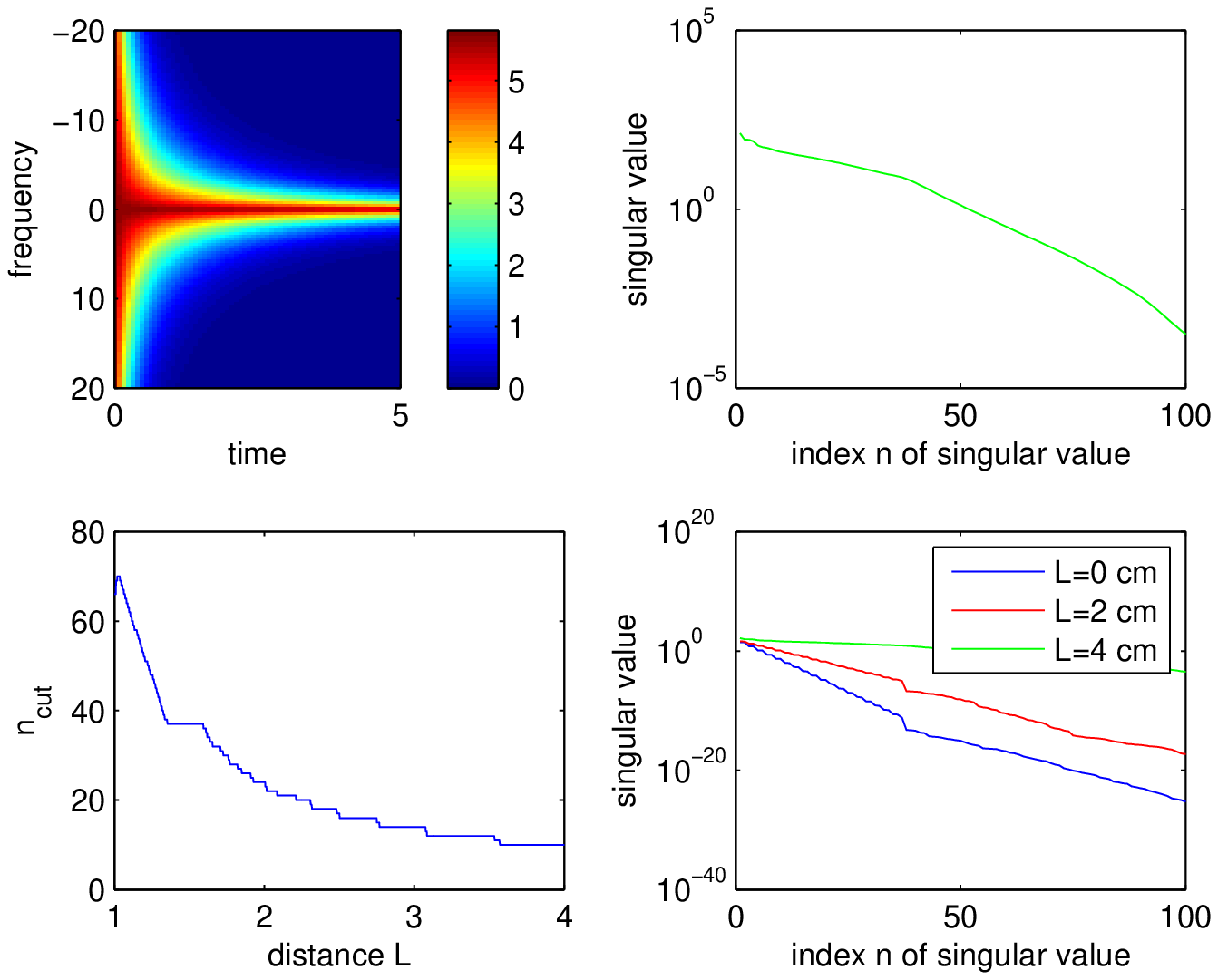}
\end{center}
\caption{Case: $\gamma=1.1$.  The upper left and right pictures visualize
the kernel $\fourier{\M}(\omega,t')$ and its singular values.
The lower right and left pictures visualize $\fourier{\M}(\omega,t')$ for the detector
distances $L=0\cdot cm$, $L=2\cdot cm$ and $L=4\cdot cm$ and the respective indices $n_{cut}$ for
which the singular values are about $0.1$ per cent of the maximal singular value.}
\label{fig:distance1}
\end{figure}
\end{example}

\begin{example}\label{exam:distance3}
An analogous numerical example as in Example~\ref{exam:distance2} for the case $\gamma=2$
is presented in Fig.~\ref{fig:distance3}. From the lower left picture of Fig.~\ref{fig:distance3},
we see that if the distance is about $2\cdot cm$ from the boundary of the object,
then only $4$ singular values are available for the numerical estimation.
\begin{figure}[htb]
\begin{center}
\includegraphics[angle=90,width=\textwidth]{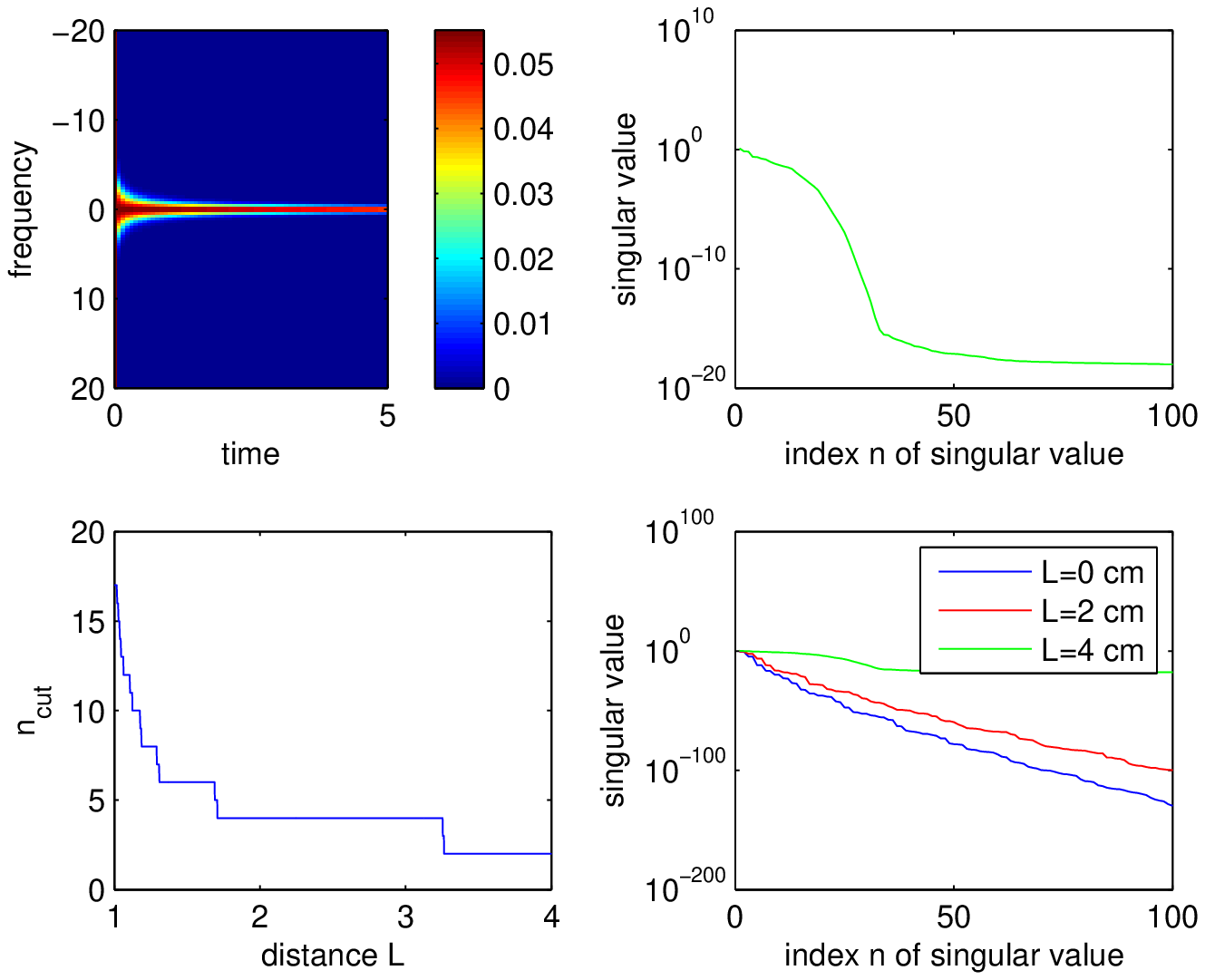}
\end{center}
\caption{Case: $\gamma=1.1$.  The upper left and right pictures visualize
the kernel $\fourier{\M}(\omega,t')$ and its singular values.
The lower right and left pictures visualize $\fourier{\M}(\omega,t')$ for the detector
distances $L=0\cdot cm$, $L=2\cdot cm$ and $L=4\cdot cm$ and the respective indices $n_{cut}$ for
which the singular values are about $0.1$ per cent of the maximal singular value.}
\label{fig:distance3}
\end{figure}
\end{example}

\begin{example}\label{exam:distance5}
An analogous numerical example as in Example~\ref{exam:distance2} for the frequency power law
$$
  \alpha_*(\omega) = \alpha_0^{pl}\cdot (-\i\,\omega)^{0.66}
  \qquad\quad (\mbox{$\alpha_0^{pl}$ as in Example~\ref{exam:distance2}})
$$
is presented in Fig.~\ref{fig:distance5}. From the lower left picture of this figure,
we see that if the distance is about $2\cdot cm$ from the boundary of the object,
then $77$ singular values are available for the numerical estimation. If the distance is about
$4cm$, then $46$  singular values are available.
\begin{figure}[htb]
\begin{center}
\includegraphics[angle=90,width=\textwidth]{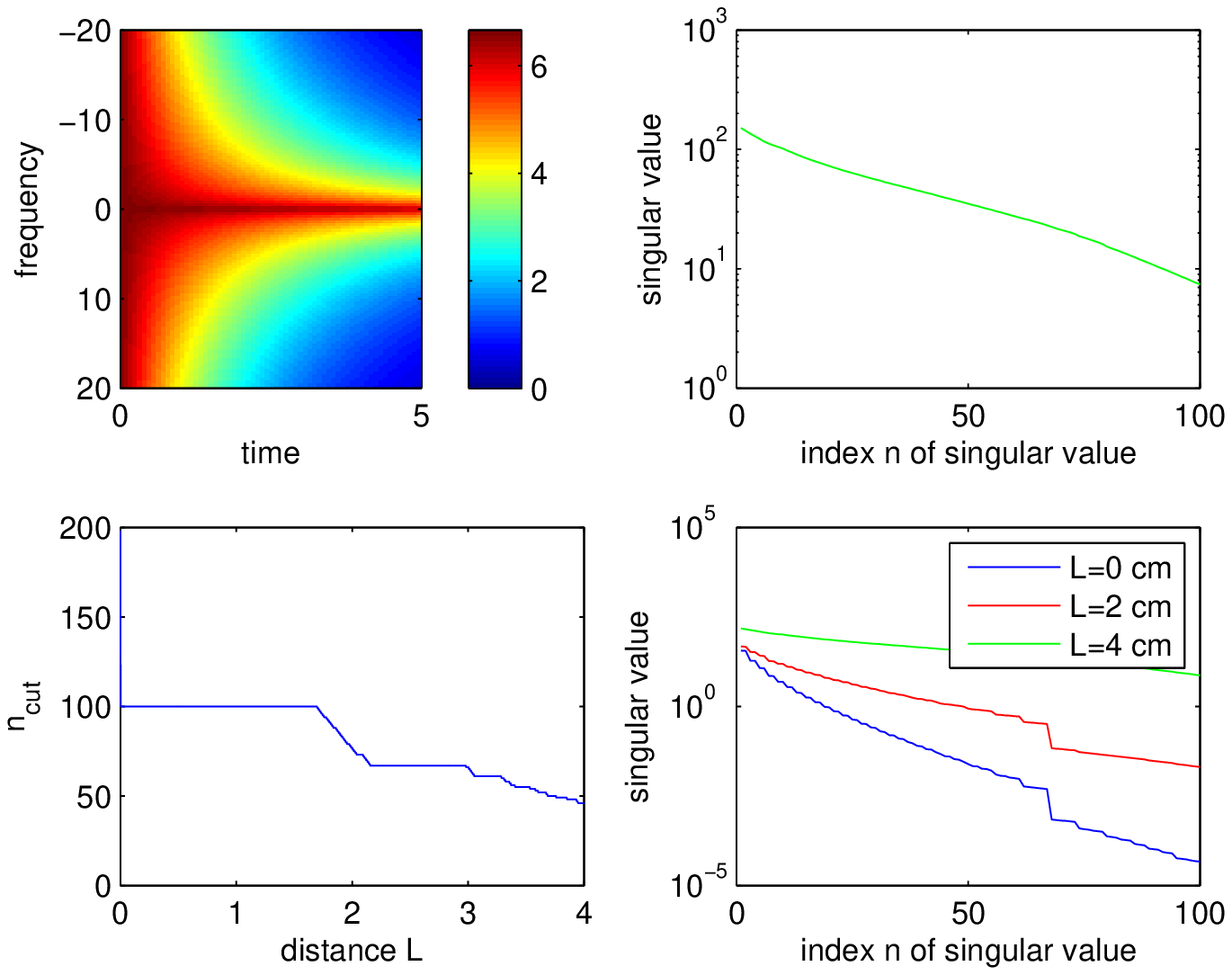}
\end{center}
\caption{Case: $\gamma=0.66$.  The upper left and right pictures visualize
the kernel $\fourier{\M}(\omega,t')$ and its singular values.
The lower right and left pictures visualize $\fourier{\M}(\omega,t')$ for the detector
distances $L=0\cdot cm$, $L=2\cdot cm$ and $L=4\cdot cm$ and the respective indices $n_{cut}$ for
which the singular values are about $0.1$ per cent of the maximal singular value.}
\label{fig:distance5}
\end{figure}
\end{example}

Comparing all numerical examples shows that the larger $\gamma$ (stronger attenuation), the more
rapidly decrease the singular values.

\section{Appendix: Nomenclature and Elementary Facts}
\label{sec:app}

\begin{description}
\item{Sets:} $B_R$ denotes the open ball with center at $\mathbf{0}$ and radius $R$. $S^n \subseteq \R^n$ denotes the $n$-dimensional unit sphere.
\item{Real and Complex Numbers:} $\C$ denotes the space of complex numbers, $\R$ the space of reals. For a complex number $c=a+\i b$
      $a=\Re{(c)}$, $b=\Im{(c)}$ denote the real and imaginary parts, respectively.
\item For a complex number $c$ we denote by $\abs{c}$ the absolute value and by $\phi \in (-\pi,\pi]$ the argument.
      That is
      \begin{equation*}
      c = \abs{c} \ef{i \phi}\;.
      \end{equation*}
      As a consequence, when $w=r\,\ef{\i \phi}$ then
      \begin{equation}\label{defpower}
        w^\gamma = \ef{\gamma\,\left(\text{log}(r)+\i \phi\right)}\;.
      \end{equation}
      Consequently $w^\gamma$ has absolute value $r^\gamma$ and the argument is $\gamma \phi$ modulo $2\pi$.
      In this paper all power functions are defined on $\C\backslash \R_-$.  We note that
      \begin{equation}\label{propsqrt}
         w\in\C\backslash \R_-   \quad\mbox{and}\quad \Re(\sqrt{w})>0  \qquad \Rightarrow \qquad
         \Im(w)\,\Im(\sqrt{w})\geq 0  \,.
      \end{equation}
\item{Differential Operators:} $\nabla$ denotes the gradient. $\nabla \cdot$ denotes divergence, and $\nabla^2$ denotes the Laplacian.
\item{Product:} When we write $\cdot$ between two functions, then it means a pointwise product, it can be a scaler product or if
the functions are vector valued an inner product. The product between a function and a number is not explicitly stated.
\item{Composition:} The composition of operators $\A$ and $\B$ is written as $\A\B$.

\item{Special functions:}
\begin{itemize}
\item The \emph{signum} function is defined by
\[
\sgn := \sgn (\vx) := \frac{\vx}{\abs{\vx}}\;.
\]
In $\R^3$ it satisfies
\begin{equation}
\label{eq:der_sgn} \nabla \cdot \sgn = \frac{2}{\abs{\vx}}\;.
\end{equation}
\item The \emph{Heaviside} function
\[
 \Heavi := \Heavi(t) := \left\{ \begin{array}{rcl}
                      0 & \text{ for } & t < 0\\
                      1 & \text{ for } & t > 0\\
                     \end{array} \right.
\]
satisfies
\[
\Heavi := \frac{1}{2} (1+\sgn)\;.
\]
\item The $\delta$-distribution is the derivative of the Heaviside function at $0$ and is denoted by
$\delta_t := \delta_t(t)$. In our terminology $\delta_t$ denotes a \emph{one}-dimensional distribution.
Sometimes, if the context is clear, we will omit the subscript at the $\delta$-distributions.
\item The three dimensional $\delta$-distribution $\delta_\vx$ is the tensor product of the three one-dimensional
distributions $\delta_{x_i}$, $i=1,2,3$. Moreover,
\begin{equation}
\label{eq:hatdelta} \delta_{\vx,t} := \delta_{\vx,t}(\vx,t) = \delta_\vx \cdot \delta_t,
\end{equation}
is a four dimensional distribution in space and time. If we do not add a subscript $\delta$ denotes a one-dimensional
$\delta$-distribution.
\item
$\chi_\Omega$ denotes the characteristic set of $\Omega$, i.e., it attains the value $1$ in $\Omega$ and is zero else.
\end{itemize}
\item{Properties related to functions:} $\supp(g)$ denote the \emph{support} of the function $g$,
that is the closure of the set of points, where $g$ does not vanish. \item{Derivative with respect to
radial components:} We use the notation
\[
 r:=r(\vx) = \abs{\vx},
\]
and denote the derivative of a function $f$, which is only dependent on the radial component
$\abs{\vx}$, with respect to $r$ (i.e., with respect to $\abs{\vx}$) by $\cdot'$.

Let $\beta = \beta(r)$, then it is also identified with the function $\beta = \beta (\abs{\vx})$
and therefore
\[
 \nabla \beta = \frac{\vx}{\abs{\vx}} \beta'\;.
\]
\item{Convolutions:} Three different types of convolutions are considered: $*_t$ and $*_\omega$
denote \emph{convolutions} with respect to time and frequency, respectively. Let $f$, $\hat{f}$,
$g$ and $\hat{g}$ be functions defined on the real line with complex values. Then
\[
\begin{aligned}
 &f *_t g := \int_\R f(t-t')g(t') d t' , \quad\quad
 &\hat{f} *_\omega \hat{g}
        := \int_\R \hat{f}(\omega-\omega')\hat{g}(\omega') d \omega'.
\end{aligned}
\]
$*_{\vx,t}$ denotes space--time convolution and is defined as follows: Let $f,g$ be functions defined
on the Euclidean space $\R^3$ with complex values, then
\[
f *_{\vx,t} g := \int_{\R^3} \int_\R f(\vx - \vx',t-t')g(\vx',t') d\vx' d t'\;.
\]
\item{Fourier transform:} For more background we refer to~\cite{Lig64,Tit48,Pap62,Yos95,Hoe03}. All along this paper $\fourier{\cdot}$
denotes the \emph{Fourier transformation} with respect to $t$, and the \emph{inverse Fourier transform} $\ifourier{\cdot}$ is with respect to $\omega$.
In this paper we use the following definitions of the transforms:
\begin{equation*}
\begin{aligned}
\fourier{f}(\omega) &= \frac{1}{\sqrt{2\pi}} \int_\R \ef{\i \omega t} f(t) d t\,,\\
\ifourier{\hat{f}}(t) &= \frac{1}{\sqrt{2\pi}} \int_\R \ef{-\i
\omega t} \hat{f}(\omega) d\omega\;.
\end{aligned}
\end{equation*}
The Fourier transform and its inverse have the following properties:
\begin{enumerate}
\item\label{item:derF}
      \[
      \fourier{\frac{\partial }{\partial t}f} (\omega) =  (-\i \omega)\fourier{f}(\omega)\;.
      \]
\item \[
      \begin{aligned}
      \fourier{f \cdot g} & =\frac{1}{\sqrt{2\pi}} \fourier{f} *_\omega \fourier{g} \text{ and }\\
      \fourier{f} \cdot \fourier{g} &= \frac{1}{\sqrt{2\pi}} \fourier{f *_t g},\\
      \ifourier{\hat{f} \cdot \hat{g}}
             &= \frac{1}{\sqrt{2\pi}}\ifourier{\hat{f}} *_t \ifourier{\hat{g}} \text{ and }\\
      \ifourier{\hat{f}} \cdot \ifourier{\hat{g}}
             &= \frac{1}{\sqrt{2\pi}} \ifourier{\hat{f} *_\omega \hat{g}}\;.
      \end{aligned}
      \]
\item   For $a \in \R$
\[
    \fourier{f(t-a)}(\omega) = \ef{i a \omega} \cdot\fourier{f(t)}(\omega)
\]
\item\label{item:Fdelta}
    The $\delta$-distribution at $a \in \R$ satisfies
    \[
   \delta_t(t-a)=\frac{1}{\sqrt{2\pi}}\ifourier{\exp(\i a \omega)}(t)\;.
\]
\item  \label{item:even}
  Let $f$ be real and even, odd respectively, then $\fourier{f}$ is real and even, imaginary and odd,
   respectively.
\item \label{item:temp} The Fourier transformation of a tempered distribution
is a tempered distribution.
\item \label{item:conv}
      Let $\tau_1,\,\tau_2\in \R$. If $f_1$ and $f_2$ are two distributions with support
      in $[\tau_1,\infty)$ and $[\tau_2,\infty)$, respectively, then $f_1*f_2$ is well-defined
      and (cf. ~\cite{Hoe03})
      \begin{equation}\label{propIM}
            \mbox{supp} (f_1 * f_2)
           \subseteq  \mbox{supp} (f_1) + \mbox{supp} (f_2) \subseteq [\tau_1+\tau_2,\infty)\,.
      \end{equation}
\end{enumerate}
\item{The Hilbert transform for $L^2-$functions is defined by}
\[
 \hilbert{f} (t) = \frac{1}{\pi} \Xint-_\R \frac{f(s)}{t-s}ds\;,
\]
where $\Xint-_\R f(s) ds$ denotes the Cauchy principal value of $\int_\R f(s) ds$.

A more general definition of the  Hilbert transform can be found in~\cite{BelWoh66}. The Hilbert
transform satisfies
\begin{itemize}
 \item $\hilbert{\fourier{f}}(\omega) = -i\fourier{\sgn f}(\omega)$,
 \item $\hilbert{\hilbert{f}} = - f$.
\end{itemize}
From the first of these properties the Kramers-Kronig relation can be formally derived as follows.
Since $f(t)$ is a causal function if and only if $ f = \Heavi \cdot f$ and $\Heavi=(1+\sgn)/2$,
it follows that $\fourier{f} =  [\fourier{f} +\i \hilbert{\fourier{f}}]/2 $, which is equivalent to
$\fourier{f} = i\hilbert{\fourier{f}}$, i.e.
\[
      \Re(\fourier{f}) =  -\Im(\hilbert{\fourier{f}})  \ltext{and}
      \Im(\fourier{f}) =  \Re(\hilbert{\fourier{f}}) .
\]
\item{The inverse Laplace transform} of $f$ is defined by
\[
\ilaplace{f}(t) = \left\{ \begin{array}{ccl}
0 &\text{ for }& t < 0,\\
\frac{1}{2 \pi i} \int_{\gamma - i\infty}^{\gamma + \i \infty} \ef{st} f(s)ds, &\text{ for }& t > 0,
\end{array}\right.
\]
where $\gamma$ is appropriately chosen.

The inverse Laplace transform satisfies (see e.g. \cite{Heu91})
\begin{equation}
\label{eq:pr1} \ilaplace{h(s-a)}(t) = \ef{a t} \ilaplace{h(s)}(t) \text{ for all } a, t \in \R
\end{equation}
and
\begin{equation}
 \label{eq:pr2}
\ilaplace{s^{-r}}(t) = \frac{H(t) t^{r-1}}{\Gamma(r)}\qquad (r>0)\;.
\end{equation}
\end{description}

\section*{Acknowledgement}
This work has been supported by the Austrian Science Fund (FWF)
within the national research network Photo\-acoustic Imaging in Biology and Medicine, project S10505-N20.

\printindex
\end{document}